\documentclass[a4paper,11pt,leqno]{article}

\usepackage{amsfonts,amsmath,amsthm,graphicx,txfonts,verbatim,yhmath,fullpage,fancyhdr}
\usepackage[all]{xy}
\xyoption{web}

\newtheorem*{thm*}{Theorem}
\newtheorem{thm}{Theorem}[section]
\newtheorem{lem}[thm]{Lemma}
\newtheorem{prop}[thm]{Proposition}
\newtheorem{cor}[thm]{Corollary}

\theoremstyle{definition}
\newtheorem{ex}[thm]{Example}
\newtheorem{rmrk}[thm]{Remark}
\newtheorem{defn}[thm]{Definition}

\newtheorem*{notation*}{Notation}
\newtheorem*{ackno*}{Acknowledgements}

\newcommand{\fS}{{\mathcal S}}
\newcommand{\triv}{\boldsymbol{1}}
\newcommand{\sign}{\epsilon}
\newcommand{\bZ}{\mathbb{Z}}
\newcommand{\bQ}{\mathbb{Q}}

\newcommand{\bC}{\mathbb{C}}
\newcommand{\bR}{\mathbb{R}}

\newcommand{\fp}{\mathfrak p}
\newcommand{\fP}{\mathfrak P}

\newcommand{\fX}{\mathfrak X}

\newcommand{\C}{{\cal C}}
\newcommand{\I}{{\cal I}}
\newcommand{\cK}{{\cal K}}
\newcommand{\cL}{{\cal L}}
\newcommand{\cR}{{\cal R}}

\newcommand{\z}{\zeta_p}

\newcommand{\ord}{{\rm ord}}

\newcommand{\Gal}{{\rm Gal}}
\newcommand{\tr}{{\rm tr}}
\newcommand{\N}{{\rm N}}
\newcommand{\Hom}{{\rm Hom}}
\newcommand{\Ext}{{\rm Ext}}
\newcommand{\sH}{{\rm H}}

\newcommand{\coker}{{\rm coker}}
\newcommand{\Norm}{{\rm Norm}}

\newcommand{\Sunits}[1]{{\eusm{O}_{S,#1}^\times}}
\newcommand{\Units}[1]{{\eusm{O}_{#1}^\times}}

\long\def\symbolfootnote[#1]#2{\begingroup
\def\thefootnote{\fnsymbol{footnote}}\footnote[#1]{#2}\endgroup}

\input cyracc.def
\newfam\cyrfam
\font\tencyr=wncyr10
\font\sevencyr=wncyr7
\font\fivecyr=wncyr5

\textfont\cyrfam=\tencyr \scriptfont\cyrfam=\sevencyr
\scriptscriptfont\cyrfam=\fivecyr

\font\teneurm=eurm10
\font\seveneurm=eurm7
\font\fiveeurm=eurm5
\newfam\eurmfam
\textfont\eurmfam=\teneurm
\scriptfont\eurmfam=\seveneurm
\scriptscriptfont\eurmfam=\fiveeurm

\font\teneusm=eusm10
\font\seveneusm=eusm7
\font\fiveeusm=eusm5
\newfam\eusmfam
\textfont\eusmfam=\teneusm
\scriptfont\eusmfam=\seveneusm
\scriptscriptfont\eusmfam=\fiveeusm
\def\eusm#1{{\fam\eusmfam\relax#1}}

\begin{document}
\title{On Brauer-Kuroda type relations of S-class numbers in dihedral extensions}
\author{Alex Bartel}

\maketitle
\begin{abstract}
Let $F/k$ be a Galois extension of number fields with dihedral Galois group of order $2q$, where
$q$ is an odd integer.
We express a certain quotient of $S$-class numbers of intermediate
fields, arising from Brauer-Kuroda relations, as a unit index. Our formula is valid
for arbitrary extensions with Galois group $D_{2q}$ and for arbitrary Galois-stable sets of primes
$S$, containing the Archimedean ones. Our results have curious applications to determining the
Galois module structure of the units modulo the roots of unity of a $D_{2q}$-extension from class
numbers and $S$-class numbers. The techniques we use are mainly representation theoretic and we
consider the representation theoretic results we obtain to be of independent interest.
\end{abstract}
\tableofcontents
\symbolfootnote[0]{MSC 2010: Primary 11R29, 11R33, Secondary 20C10, 11R27}
\newpage
\pagestyle{fancy}
\fancyhf{} 
\fancyhead[C]{\textsc{Alex Bartel -- Relations between class numbers in dihedral extensions}}
\fancyfoot[C]{\thepage}
\section{Introduction}
Dirichlet \cite{Dir-42} was the first to establish a relation between class numbers of a number field
and its subfields in 1842: he showed that for a positive integer $d$ that is not a square,
the quotient of the class number $h$ of $\bQ(\sqrt{d},\sqrt{-1})$ by the product of the class numbers
$h_d$ of $\bQ(\sqrt{d})$ and $h_{-d}$ of $\bQ(\sqrt{-d})$ is either equal to $1$ or $2$. In 1950,
Brauer \cite{Bra-51} and Kuroda \cite{Kur-50} independently initiated a systematic study of relations
between class numbers in number fields arising from isomorphisms of permutation representations of
finite groups. More precisely, if $G$ is a finite group and $\{H_i\}_i$ and $\{H_j'\}_j$ are sets of
subgroups such that there is an isomorphism of permutation representations of $G$
\[\bigoplus_i \bQ[G/H_i] \cong \bigoplus_j \bQ[G/H_j'],\]
and if $F/K$ is a Galois extension of number fields with Galois group $G$, then Artin formalism for
Artin $L$-functions implies that we have an equality of zeta-functions of the corresponding fixed fields:
\[\prod_i \zeta_{F^{H_i}}(s) = \prod_j \zeta_{F^{H_j'}}(s).\]
More generally, if $S$ is any finite $G$-stable set of places of $F$ containing all the Archimedean ones,
then we have an analogous equality of $S$-zeta functions.
Invoking the analytic class number formula (see e.g. \cite[Chap. I, Cor. 2.2]{Tat-84}) yields the
equality
\begin{eqnarray}\label{eq:classNoFormula}
\prod_i \frac{h_S(F^{H_i})R_S(F^{F_i})}{w(F^{H_i})} = \prod_j \frac{h_S(F^{H_j'})R_S(F^{F_j'})}{w(F^{H_j'})},
\end{eqnarray}
where for a number field $M$, $h_S(M)$, $R_S(M)$, and $w(M)$ denotes the $S$-class number of $M$, the
$S$-regulator of $M$, and the number of roots of unity in $M$, respectively. See below for precise
definitions.

Sometimes, the value of the class number quotient
can be given an interpretation in terms of a unit index. In Dirichlet's case, the quantity
$2h/(h_dh_{-d})$ is the index in the unit group of $\bQ(\sqrt{d},\sqrt{-1})$ of the subgroup generated by
the roots of unity and the unit group of $\bQ(\sqrt{d})$. If one wants to make an analogous statement for
a general base field and any bi-quadratic extension, then the class number quotient must have
the class numbers of
all three intermediate quadratic extensions in the denominator and the formula is more complicated due a
larger unit rank. A correct formula for bi-quadratic extensions in this more general case, with $S$
equal to the set of Archimedean primes, was only given in 1994 by Lemmermeyer \cite{Lem-94}.

Our main result is a unit index formula for Galois extensions with Galois group $D_{2q}$
for $q$ any odd integer. Let $\Units{M}$ denote the units in the ring of integers of a number field $M$.
In 1977, Halter-Koch showed:
\begin{thm*}[\cite{HK-77}, section 4]
Let $F/\bQ$ be a Galois extension with Galois group $D_{2p}$ for $p$ an odd prime. Let $K$ be the
quadratic subfield and $L\neq L'$ two intermediate extensions of degree $p$ over $\bQ$. Let $r(K)$ be the rank
of the units in $K$, which is either 0 or 1. Then
\[\frac{h(F)p^{r(K)+1}}{h(K)h(L)^2}=[\Units{F}:\Units{L}\Units{L'}\Units{K}].\]
\end{thm*}
This was generalised to arbitrary base fields by Lemmermeyer in 2005 under a
restrictive assumption on the extension:
\begin{thm*}[\cite{Lem-05}, Theorem 2.2]
Let $F/k$ be a Galois extension of number fields with Galois group $D_{2p}$ for $p$ an odd prime,
let $K$ be the intermediate quadratic extension and $L\neq L'$ two intermediate extensions of
degree $p$ over $k$. Assume that $F/K$ is unramified.
Let $r(k)$ and $r(K)$ denote the ranks of the unit groups in the respective fields. Then
\[\frac{h(F)h(k)^2p^{r(K)+1-r(k)}}{h(K)h(L)^2}=[\Units{F}:\Units{L}\Units{L'}\Units{K}].\]
\end{thm*}
In 2008, Caputo and Nuccio derived a formula for $D_{2q}$ extensions where $q$ is any odd integer
for certain base fields and certain extension:
\begin{thm*}[\cite{CN}, Theorem 3.4]
Let $k$ be a totally real number field, $F/k$ a totally imaginary Galois extension with Galois group
$D_{2q}$ where $q$ is an odd integer. Let $K$ be the intermediate quadratic and $L,L'$ fixed fields
of elements $\sigma,\sigma'$ of order 2 such that $\sigma\sigma'^{-1}$ is of order $q$. Then
\[\frac{h(F)h(k)^2q}{h(K)h(L)^2} = \frac{[\Units{F}:\Units{L}\Units{L'}\Units{K}]}{[\Units{L}\Units{L'}\cap \Units{K}:\Units{k}]}.\]
\end{thm*}

In this paper we complete the study of unit index formulae for dihedral extensions of degree $2q$.
We will only state the formula explicitly for $D_{2p}$, where $p$ is a prime and explain how it is derived for $D_{2q}$ for arbitrary odd integers $q$, since
the formula gets unwieldy in the general case, although conceptually not difficult.
\begin{thm}\label{thm:unitIndex}
Let $F/k$ be a Galois extension of number fields with Galois group $D_{2p}$ for $p$ an odd prime,
let $K$ be the intermediate quadratic extension and $L,L'$ distinct intermediate extensions of degree $p$.
Let $S$ be a finite $\Gal(F/k)$-stable set of places of $F$ including the Archimedean ones. We write $\eusm{O}_S^\times$ for
$S$-units, $h_S$ for $S$-class
numbers and $r_S$ for the ranks of $S$-units. Let $a(F/k,S)$ be the number of primes of $k$ which
lie below those in $S$ and whose decomposition group is equal to $D_{2p}$. Finally, set $\delta$ to be
3 if $L/k$ is obtained by adjoining the $p$-th root of a non-torsion $S$-unit (thus so is $F/K$)
and 1 otherwise. Then we have
\begin{eqnarray*}
\frac{h_S(F)h_S(k)^2}{h_S(K)h_S(L)^2} = p^{\alpha/2}\times
[\Sunits{F}:\Sunits{L}\Sunits{L'}\Sunits{K}],
\end{eqnarray*}
where $\alpha = 2r_S(k)-r_S(K)-\frac{r_S(F)-r_S(K)}{p-1}+a(F/k,S)-\delta$.
\end{thm}
Note that all the terms in the exponent of $p$ are very easy to compute in practice (e.g. taking $S$
to be the set of Archimedean places forces $a(F/k,S)=0$; see section \ref{sec:examples} for more examples).

For arbitrary Galois extensions and sets of subgroups $H_i$, $H_j'$ giving isomorphic permutation
representations, Brauer showed that the class number quotient
$\prod_i h(F^{H_i})/\prod_j h(F^{H_j'})$ takes only finitely many values as $F$ ranges over all Galois
extensions of $K$ with Galois group $G$ (see \cite[Satz 5]{Bra-51}). He further showed that $\prod_i w(F^{H_i})/\prod_j w(F^{H_j'})$
is a power of $2$ (see \cite[\S 2]{Bra-51}) and observed that if $p$ is a prime number not dividing the order of $G$, then
$\ord_p\left(\prod_i R(F^{H_i})/\prod_j R(F^{H_j'})\right) = 0$ (\cite[Satz 4, Bemerkung 2]{Bra-51}). However, there is to date no general
formula that explains exactly, what values the regulator quotient can take.
As a by-product of our calculations, we get the following result in this direction:
\begin{thm}\label{thm:noP}
Let $G$ be a finite group, let $N$ be a normal subgroup such that $G/N$ is cyclic, let $l$ be a prime
number not dividing the order of $N$. Let $F/K$ be a Galois extension of number fields with Galois
group $G$ and $\{H_i\}$, $\{H_j'\}$ be sets of subgroups yielding an isomorphism of permutation
representations as above. Let $S$ be a finite $G$-stable set of places of $F$ including all the Archimedean
ones. Then
\[\ord_l\left(\prod_i R_S(F^{H_i})\big/\prod_j R_S(F^{H_j'})\right) = 0,\]
where $R_S$ denotes regulators of $S$-units,
and we have an equality
of the $l$-parts of $S$-class numbers: $\prod_i h_S(F^{H_i})_l = \prod_j h_S(F^{H_j'})_l$.
\end{thm}
We will briefly describe the structure of the paper and the main ideas of the proofs.

Already Brauer pointed out that the regulator quotient is a purely representation theoretic invariant
of the $\bZ[G]$-module $\Sunits{F}$. This observation was crucial for proving that
the regulator quotient takes only finitely many values for a fixed base field and varying Galois
extensions $F$ with Galois group $G$. The main step towards the proof of both Theorem \ref{thm:unitIndex}
and Theorem \ref{thm:noP} is a representation theoretic description of the regulator quotient. We will provide
such a description in section \ref{sec:regConsts} by linking the regulator quotients to certain invariants,
first introduced by Tim and Vladimir Dokchitser in \cite{TVD-1} and further explored by the Dokchitser
brothers in \cite{TVD-2} and by the author in \cite{Bar-10} in the context of elliptic curves.
These invariants are rational numbers that can be attached to pairs consisting of an integral
representation of a group and an isomorphism of permutation representations. We will call these
numbers Dokchitser constants (deviating from the original name 'regulator constants'). To express the
regulator quotients in terms of Dokchitser constants is not entirely trivial and is done in
Proposition \ref{prop:newReg}.

De Smit \cite[Theorem 2.2]{Smi-01} has derived a different expression for the regulator quotient, 
which turns out to be closely related to ours.
In section \ref{sec:altDefn} we will give an alternative definition of the Dokchitser constants
and will show how this ties in with de Smit's result. The alternative definition will also be useful
to derive some properties of the Dokchitser constants, which will lead to a proof of Theorem
\ref{thm:noP}.
We think that this section is of independent interest, since it sheds some light on the nature of
Dokchitser constants and therefore complements the results of \cite{TVD-2}.

In section \ref{sec:dihGroups} we turn to Dokchitser constants in dihedral groups. As it turns out,
one can compute all the Dokchitser constants for all integral representations of $D_{2p}$. We should
mention that $D_{2p}$ must be regarded as a lucky exception. Although it suffices
to determine the Dokchitser constants for indecomposable representations (see Proposition
\ref{prop:multiplicativity}), a finite group can have
infinitely many non-isomorphic indecomposable integral representations and nobody knows how to
classify them in general. However, $D_{2p}$ has finitely many and they have been written down
explicitly in \cite{Lee-64}. Still, it is not clear a priori that their Dokchitser constants
can be computed in general, since their number grows with $p$.

In section \ref{sec:mainRes} we use the properties of Dokchitser constants established in section
\ref{sec:altDefn} to prove Theorem \ref{thm:noP}.
We then use the computation of Dokchitser constants for $D_{2p}$ to prove
Theorem \ref{thm:unitIndex} and explain how, using formal properties of Dokchitser constants,
one can derive a formula for $D_{2q}$-extensions for any odd number $q$. Surprisingly, the generalisation to $D_{2q}$
is rather easy. Because the most general formula would look rather long and obstruct its conceptual
simplicity, we will not write it down. Considering $S$-units instead of just units also introduces
very few conceptual difficulties. We note that the way we obtain a unit index formula for
$D_{2q}$-extensions is a completely general procedure to glue unit index formulae together from
intermediate extensions. We hope that this will prove useful in the search of unit index formulae
in much more general contexts.

In the last section we give various examples. For example, we show how the formulae of
Halter-Koch and of Lemmermeyer follow from our Theorem \ref{thm:unitIndex}. We also demonstrate how
our computations can sometimes be used to determine the structure of the Galois module given by the
units modulo torsion in a dihedral extension in terms of the class numbers and $S$-class numbers
of the field and its subfields.

We should mention that we use very little number theory in this paper. We need the analytic class
number formula (or merely its compatibility with Artin formalism), but unlike the proof of a special
case of Theorem \ref{thm:unitIndex} in \cite{Lem-05}, we do not need any class field theory.
\begin{ackno*}
 This research was done, while I was a member of the Department of
Pure Mathematics and Mathematical
Statistics at the University of Cambridge and of St. John’s College, Cambridge.
I would like to thank both institutions for support of very diverse nature! Many thanks are due to
Tim and Vladimir Dokchitser and to Antonio Lei for many helpful discussions and to Vladimir for many
helpful comments. I am also very grateful to
Samir Siksek for an idea which proved crucial in giving Definition \ref{defn:pairing} and to Luco
Caputo for pointing out two mistakes in an earlier version of the manuscript. Finally, thanks a due to
an anonymous referee for carefully reading an earlier version and for many useful remarks and
comments. I gratefully acknowledge the financial support through an EPSRC grant.
\end{ackno*}
\begin{notation*} Throughout the paper, the following notation will be used for a number field $F$
and for $S$ a finite set of places of $F$ containing the Archimedean ones:
\newline
\begin{tabular}{ l l }
  $h_S(F)$ & the $S$-class number of $F$, i.e. the class number of the ring of all\\
  & elements of $F$ which are integral at all places not in $S$.\\
  $w(F)$ & the number of roots of unity in $F$.\\
  $\Sunits{F}$ & the group of $S$-units of $F$.\\
  $r_S(F)$ & the $\bZ$-rank of $\Sunits{F}$, i.e. $|S|-1$.\\
  $U_S(F)$ & the group of $S$-units modulo torsion; we will often identify units of $F$\\
  & with their image in $U_S(F)$, when no confusion can arise.\\
  $R_S(F)$ & the $S$-regulator of $F$, i.e. the absolute value of the determinant of the \\
  & square matrix
   of size $r_S(F)$, whose $(i,j)$-th entry is $||\log(u_i)||_{\fp_j}$ where\\
  & $\fp_j$ runs through
  the set of all but one absolute values attached to the \\
  & places in $S$ and
  $\{u_1,\ldots,u_{r_S(F)}\}$ is a set of generators of the group of\\
  & $S$-units mod torsion.\\
  $\zeta_{F,S}(s)$ & the $S$-zeta function of $F$, $\zeta_{F,S}(s)= \prod_{\fp\notin S} (1-\N\fp^{-s})^{-1}$\\
  & for $\Re(s) > 1$, the product taken over the places of $F$ not in $S$ and\\
  & $\N\fp$ denoting the absolute norm of $\fp$.\\
\end{tabular}
\newline
The normalisations of the absolute values $||.||_\fp$ attached to places $\fp$ are as follows:
if $F_\fp=\bR$,
then the absolute value is just the usual real absolute value. If $F_\fp=\bC$, then it is the square of
the usual absolute value. If $F_\fp$ is a $p$-adic field with residue field of size $q$, then the
attached
absolute value is the $p$-adic absolute value, normalised in such a way that $||\pi||_\fp=1/q$ for any
uniformiser $\pi$ in $F_\fp$.

When $k$ is a subfield of $F$, we will write $S|_k$ for the set of places of $k$ lying below those in
$S$. We will often write $h_S(k)$ etc. instead of $h_{S|_k}(k)$.
\end{notation*}

\section{Regulator quotients and Dokchitser constants}\label{sec:regConsts}
In this section we recall the definition of Dokchitser constants from \cite{TVD-1} (where they were
called regulator constants) and relate them to quotients of regulators of number fields. But first,
we introduce a convenient language to talk about the Brauer-Kuroda type relations.

\subsection{Relations of permutation representations and Dokchitser constants}
Let $G$ be any finite group. We recall the following standard definitions (see e.g. \cite{CR}):
\begin{defn}
The \emph{Burnside ring} of $G$ is defined as the ring of formal $\bZ$-linear combinations of isomorphism
classes $[X]$ of finite $G$-sets modulo the relations
$$[X] + [Y] = [X\sqcup Y],\;\;\;\;[X][Y] = [X\times Y],$$
where $X\sqcup Y$ denotes the disjoint union and $X\times Y$ denotes the Cartesian product.
\end{defn}
The set of isomorphism classes of transitive $G$-sets is in bijection with the set of conjugacy
classes of subgroups of $G$ via the map which assigns to the subgroup $H$ the set of co-sets $G/H$.
We will usually represent elements of the Burnside ring as formal sums $\sum_i H_i - \sum_j H_j'$ using
this identification.
\begin{defn}
Let $A$ be either $\bQ$ or $\bZ_{(p)}$, the localisation of $\bZ$ at a prime $p$. The \emph{representation
ring} of $G$ over $A$ is the ring of formal $\bZ$-linear combinations of
isomorphism classes $[M]$ of $A$-free finite dimensional $AG$-modules modulo the relations
$$[M] + [N] = [M\oplus N],\;\;\;\;[M][N] = [M\otimes N].$$
Here and in the rest of the paper, $AG$ denotes the group algebra of the group $G$ over $A$.
\end{defn}
We have a natural map from the Burnside ring to the representation ring that sends a $G$-set $X$ to
the $AG$-module $A\left[X\right]$ with $A$-basis indexed by the elements of $X$ and the natural
$G$-action. If we take $A$ to be $\bQ$, then the image of the Burnside ring in the representation
ring has finite index (called the Artin index of the group $G$).
\begin{defn}\label{defn:relation}
We will call an element $\Theta$ of the kernel of the above map from the Burnside ring of $G$ to the
representation ring over $A$ an \emph{$AG$-relation}. If $A=\bQ$, then we will drop $A$ from the notation
and will just say that $\Theta$ is a $G$-relation.
\end{defn}
The number of isomorphism classes of irreducible rational representations of a finite group $G$ is equal
to the number of conjugacy classes of cyclic subgroups of $G$ (see \cite[\S 13.1, Cor. 1]{Ser-67}).
Since, as remarked above, the image of the Burnside ring has full rank in the representation ring over $\bQ$, 
the lattice of $G$-relations has rank equal to the number of conjugacy classes of non-cyclic
subgroups.
\begin{ex}\label{ex:d2p}
Let $p$ be an odd prime. The dihedral group $D_{2p}$ with $2p$ elements has one non-cyclic subgroup,
namely itself. Decomposing the various permutation representations into irreducible summands, one easily
finds that $\Theta = 1 - 2C_2 - C_p + 2D_{2p}$ is a relation. Since it is not divisible by any
integer, it must span the $\bZ$-lattice of $D_{2p}$-relations.
\end{ex}
We now recall the concept that will be of central importance in this paper:
\begin{defn}\label{defn:dokConst}
Let $G$ be a finite group, let $\Theta=\sum_i H_i - \sum_j H_j'$ be an $AG$-relation and let $\cR$ be
a principal ideal domain such that its field of fractions $\cK$ has characteristic not dividing $|G|$.
Given an $\cR$-free finite rank $\cR G$-module $\Gamma$ such that $\Gamma\otimes \cK$ is self-dual
we fix a non-degenerate $G$-invariant bilinear pairing $\left<,\right>$ on $\Gamma$ with values in some
extension $\cL$ of $\cK$. For any subgroup
$H$ of $G$, the fixed points $\Gamma^H$ are also $\cR$-free, since $\cR$ is a PID, and the pairing is
also non-degenerate when restricted to $\Gamma^H$ by \cite[Lemma 2.15]{TVD-2}. We may thus define
the \emph{Dokchitser constant} of $\Gamma$ with respect to $\Theta$ to be
$$\C_\Theta(\Gamma) = \frac{\prod_i \text{det}\left(\frac{1}{|H_i|}\left<,\right>\big|\Gamma^{H_i}\right)}
{\prod_j \text{det}\left(\frac{1}{|H'_j|}\left<,\right>\big|\Gamma^{H'_j}\right)} \in \cL^\times/(\cR^\times)^2,$$
where each inner product matrix is evaluated with respect to some $\cR$-basis on the fixed submodule.
If the matrix of the pairing on $\Gamma^{H}$ with respect to some fixed basis is $M$, then changing
the basis by the change of basis matrix $X\in \text{GL}(\Gamma^{H})$ changes the matrix of the pairing
to $X^{\rm tr}MX$. So the Dokchitser constant is indeed a well-defined element of $\cL^\times/(\cR^\times)^2$.
\end{defn}
\textbf{Convention.} From now on, $\cR$ will be assumed to be a PID with field of fractions $\cK$ of
characteristic not dividing $|G|$. We will refer to $\cR G$-modules $\Gamma$ that are free and of finite
rank over $\cR$ as $\cR G$-lattices. We will always assume that $\Gamma\otimes \cK$ is self-dual.
When we refer to
subgroups we will always mean subgroups up to conjugation, unless specifically otherwise stated.
So the subgroups $H$ and $H'$ will be treated as being the same if the $G$-sets $G/H$ and $G/H'$ give
the same element of the Burnside ring.

The choice of pairing is not present in the notation of
Dokchitser constants and indeed we have:
\begin{thm}[\cite{TVD-2}, Theorem 2.17]\label{thm:indepPairing}
The value of ${\C}_\Theta(\Gamma)$ is independent of the choice of the pairing.
\end{thm}
In particular, the pairing can always be chosen to be $\cK$-valued and so we see that the Dokchitser
constant is in fact an element of $\cK^\times/(\cR^\times)^2$. Note that if $\cR=\bZ$, then the
Dokchitser constant is just a rational number. If $\cR=\bZ_p$, then at least the $p$-adic order of
the Dokchitser constant is well-defined. If on the other hand $\cR=\bQ$, then the Dokchitser constant
is only defined up to rational squares, and if $\cR=\bQ_p$, then only the parity of the $p$-adic order
is defined.
An immediate consequence of Theorem \ref{thm:indepPairing} is
\begin{prop}[\cite{TVD-2}, Corollary 2.18]\label{prop:multiplicativity}
The Dokchitser constants are multiplicative in $\Theta$ and in $\Gamma$, i.e.
\begin{eqnarray*}{\C}_\Theta(\Gamma\oplus\Gamma') = \C_\Theta(\Gamma)\C_\Theta(\Gamma'),\\
\C_{\Theta+\Theta'}(\Gamma) = \C_\Theta(\Gamma)\C_{\Theta'}(\Gamma).
\end{eqnarray*}
\end{prop}
\begin{ex}
Take $G=S_3$. There are three irreducible complex representations of $S_3$, namely the trivial
representation $\triv$, the one-dimensional sign representation $\sign$ and a two-dimensional
representation $\rho$, and they are all defined
over $\bQ$. We saw in Example \ref{ex:d2p} that there is, up to integer multiples, a unique relation
$$1 - 2C_2 - C_3 + 2S_3$$
and it is easy to check that the corresponding Dokchitser constants (with $\cR=\bQ$) of all three irreducible representations
are equal to 3 modulo rational squares. The representations $\triv$ and $\epsilon$ contain a unique
$G$-invariant $\bZ$-lattice each up to isomorphism and their Dokchitser constants (with $\cR=\bZ$) are $1/3$ and $3$,
respectively. The 2-dimensional representation $\rho$ contains two non-isomorphic $G$-invariant 
$\bZ$-lattices. Both can be visualised as hexagonal lattices, generated by two shortest distance
vectors $P$ and $Q$, on which the 3-cycles act as rotations by $120^{\circ}$.
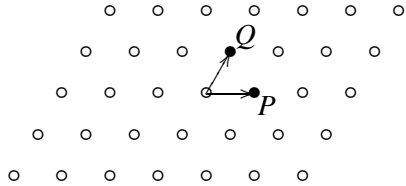
\[
\def\latticebody{%
\ifnum\latticeA=1 \ifnum\latticeB=0%
\else\drop{\circ}\fi\else
\ifnum\latticeA=0 \ifnum\latticeB=1%
\else\drop{\circ}\fi\else\drop{\circ}\fi\fi}
\xy *\xybox{0;<1.5pc,0mm>:<0.75pc,1.299pc>::
,0,{\xylattice{-3}3{-2}2}
 ,(0,0)="O"\ar
 ,(1,0)="P"*{\bullet}*+<2pt>!UL{P}\ar
 ,(0,1)="Q"*{\bullet}*+<3pt>!DL{Q}
}
\endxy
\]
On one, the 2-cycles act by reflection through a shortest distance vector (eg. through $P$) and on
the other the 2-cycles act by reflection through the long diagonal of the fundamental parallelograms
(which are $P+Q$ and its rotations by $120^{\circ}$ in the sketch). Each one of the two can be
embedded into the other $G$-equivariantly with index 3, but there is no $G$-equivariant bijection
between them. The Dokchitser constants (again with $\cR=\bZ$) of the two lattices are easily computed to be 1/3 and 3, respectively.
\end{ex}
\subsection{Some properties of Dokchitser constants}
We will collect some properties of Dokchitser constants that we will need later. The details can
be found in \cite{TVD-2}.
We first quote a result that shows
that, at least for $\bQ G$-modules, only \emph{finitely many} primes $p$ can appear in the Dokchitser
constants:
\begin{prop}\label{prop:notDivble}
If $\cR=\bQ$ or $\bQ_p$ and $p\nmid |G|$, then $\ord_p(\C_\Theta(\Gamma))$ is even for any
$G$-relation $\Theta$.
\end{prop}
\begin{proof}
See \cite[Corollary 2.28]{TVD-2}.
\end{proof}
In section \ref{sec:altDefn} we will generalise this statement to $\cR=\bZ$ and $\cR=\bZ_p$ and we
will further restrict for what primes $\ord_p(\C_\Theta(\Gamma))$ can be non-zero.

Relations can be restricted to subgroups, induced from subgroups and lifted from quotients as follows: 
let $\Theta=\sum_i H_i - \sum_j H_j'$ be a $G$-relation.
\begin{itemize}
\item \textbf{Induction.} If $G'$ is a group containing $G$, then by transitivity of induction,
$\Theta$ can be induced to a $G'$-relation $\Theta\!\!\uparrow^{G'}=\sum_i H_i - \sum_j H_j'$ of $G'$.
\item \textbf{Inflation.} If $G\cong \tilde{G}/N$, then each $H_i$ corresponds to a subgroup $\tilde{H_i}$
of $\tilde{G}$ containing $N$ and similarly for $H_j'$ and, inflating the permutation representations
from a quotient, we see that $\tilde{\Theta}=\sum_i \tilde{H_i} -\sum_j \tilde{H_j'}$
is a $\tilde{G}$-relation.
\item \textbf{Restriction.} If $H$ is a subgroup of $G$, then by Mackey decomposition $\Theta$ can be
restricted to an $H$-relation $\Theta\!\!\downarrow_H=\displaystyle{\sum_i \sum_{g\in H_i\backslash G/H}} H\cap H_i^g -
\displaystyle{\sum_j \sum_{g\in H_j'\backslash G/H}} H\cap H_j'^g$.
\end{itemize}
We have the following compatibility between these operations and the corresponding operations applied
to representations $\Gamma$:
\begin{prop}\label{prop:restrictionInduction}
Let $G$ be a finite group and $\Gamma$ an $\cR G$-lattice.
\begin{itemize}
\item If $H<G$ and $\Theta$ is an $H$-relation, then $\C_\Theta(\Gamma\!\!\downarrow_H) =
\C_{\Theta\uparrow^G}(\Gamma)$
\item If $G\cong \tilde{G}/N$ and $\Theta$ is a $G$-relation with $\tilde{\Theta}$ the lifted relation,
then $\C_\Theta(\Gamma) = \C_{\tilde{\Theta}}(\Gamma)$ where $\Gamma$ can also be regarded as a
$\tilde{G}$-representation.
\item If $G<G'$ and $\Theta$ is a $G'$-relation, then $\C_\Theta(\Gamma\!\!\uparrow^{G'}) = 
\C_{\Theta\downarrow_G}(\Gamma)$.
\end{itemize}
\end{prop}
\begin{proof} See \cite[Proposition 2.45]{TVD-2}.
\end{proof}
\subsection{Quotients of regulators of number fields}\label{sec:newReg}
We now want to explain the relationship between Dokchitser constants and quotients of regulators in
Brauer-Kuroda type relations. Let $G$ be a finite group, $\Theta = \sum_i H_i - \sum_j H_j'$ a
$G$-relation and $F/k$ a Galois extension of number fields with Galois group $G$. Let $S$ be a finite
$G$-stable set of places of $F$ including all the Archimedean ones.

In the definition of the Dokchitser constant we need to fix a pairing on our $\bZ G$-lattice, so to
turn regulator quotients into Dokchitser constants, we need to find a suitable pairing on $U_S(F)$. It seems
tempting therefore to multiply the matrix $(\log||u_i||_{\fp_j})$, whose determinant is the $S$-regulator of
the field, with its transpose and to take the pairing of which the resulting matrix will be the Gram
matrix. In other words, we would then define the inner product of $u_i$ with $u_j$ by
$\sum_{k=1}^{|S|-1} (\log||u_i||_{\fp_k}\log||u_j||_{\fp_k})$, with the sum running over all but one place in $S$.
This approach does not work because the resulting pairing does not behave well upon restriction to
subfields. We need its restriction to a subfield $M$ to be equal to the pairing of that subfield scaled
by the degree of $F/M$. We are very grateful to Samir Siksek for suggesting to us to try instead
summing over all places in $S$, rather than all but one. We also need another slight twist:
\begin{defn}\label{defn:pairing}
Let $M$ be a number field and $\fS$ a finite set of places of $M$ including the Archimedean ones.
Define the bilinear pairing $\left<,\right>_M$ on $U_\fS(M)$ by
\[\left<u_i,u_j\right>_M = \sum_{\fP\in \fS}\frac{1}{e_\fP f_\fP}\log||u_i||_\fP\log||u_j||_\fP,\]
where $e_\fP$ is the absolute ramification index of $\fP$ and $f_\fP$ is the degree of its residue
field over the prime subfield (defined to be 1 if $\fP$ is Archimedean).
\end{defn}
We begin by establishing the non-degeneracy of the pairing and by linking its determinant to the usual
$\fS$-regulator of a number field.
\begin{lem}\label{lem:detPairing}
Let $M$ be a number field and $\fS$ a finite set of places of $M$ containing all the Archimedean ones.
Then we have
\[\det\left(\left<,\right>_M|U_\fS(M)\right) =
\frac{\sum_{\fP\in \fS}f_\fP e_\fP}{\prod_{\fP\in \fS}f_\fP e_\fP}\cdot R_\fS(M)^2.\]
\end{lem}
\begin{proof}
Write $\fS=\{\fP_1,\ldots,\fP_{r+1}\}$ and define the following matrix:
\[\fX=\left(
\begin{array}{c c c c}
\sqrt{f_{\fP_1} e_{\fP_1}} & \frac{\displaystyle \log||u_1||_{\fP_1}}{\sqrt{\displaystyle f_{\fP_1}e_{\fP_1}}}
  & \cdots & \frac{\displaystyle \log||u_r||_{\fP_1}}{\sqrt{\displaystyle f_{\fP_1}e_{\fP_1}}}\\
\vdots & \vdots & \ddots & \vdots\\
\sqrt{f_{\fP_{r+1}} e_{\fP_{r+1}}} & \frac{\displaystyle \log||u_1||_{\fP_{r+1}}}{\sqrt{\displaystyle f_{\fP_{r+1}}e_{\fP_{r+1}}}}
  & \cdots & \frac{\displaystyle \log||u_r||_{\fP_{r+1}}}{\sqrt{\displaystyle f_{\fP_{r+1}}e_{\fP_{r+1}}}}\\
\end{array}
\right).\]
Then by the product formula and by the definition of our pairing, we have
\[\fX^\tr \fX=\left(
\begin{array}{c c c c}
\sum_{\fP\in \fS}f_\fP e_\fP & 0 & \cdots & 0\\
0 & \left<u_1,u_1\right>_M & \cdots & \left<u_1,u_r\right>_M\\
\vdots & \vdots & \ddots & \vdots\\
0 & \left<u_r,u_1\right>_M & \cdots & \left<u_r,u_r\right>_M\\
\end{array}
\right)\]
and so 
\begin{eqnarray}\label{eq:AtrA}
\det(\fX^\tr \fX)= \det\left(\left<,\right>_M|U_\fS(M)\right)\cdot \sum_{\fP\in \fS}f_\fP e_\fP.
\end{eqnarray}
On the other
hand, thanks to the product formula, by multiplying the $i$-th row of $\fX$ by $\sqrt{f_{\fP_i}e_{\fP_i}}$ for each $i$ and by adding
all the rows of the resulting matrix to the last one, we get zeros in the entire bottom row, apart from
the bottom left entry, where we get $\sum_{\fP\in \fS} f_{\fP}e_{\fP}$. Moreover, the resulting matrix
with the first column and the bottom row deleted has determinant equal to $R_\fS(M)$. In summary we see
that
\begin{eqnarray}\label{eq:A}
\det(\fX) = \frac{\sum_{\fP\in \fS}f_{\fP}e_{\fP}}{\prod_{\fP\in \fS} \sqrt{f_{\fP}e_{\fP}}} \cdot R_\fS(M),
\end{eqnarray}
and by combining equations (\ref{eq:AtrA}) and (\ref{eq:A}), the claim follows.
\end{proof}
We now return to our previous scenario and explain how the pairing behaves in relations. It is clear
that if $F/k$ is a Galois extension of number fields and $S$ is a finite set of places of $F$ containing
all the Archimedean ones, then $\left<,\right>_F$ is $G$-invariant. It also behaves correctly under
restriction to subfields:
\begin{lem}\label{lem:baseChange}
Let $F/k$ be a Galois extension, $S$ a finite Galois stable set of places of $F$ including the Archimedean
ones, $L\leq M$ subfields of $F$ containing $k$ and $u_i,u_j$ $S$-units in $U_S(L)$. Then
\[\left<u_i,u_j\right>_M = \left[M:L\right]\left<u_i,u_j\right>_L.\]
\end{lem}
\begin{proof}
This is easy to see by considering each prime of $S|_L$ separately, since, with our normalisations
of the absolute values, we have
\[\log||u||_\fP = e_{\fP/\fp}f_{\fP/\fp}\log||u||_\fp\]
for any $u\in U_S(L)$ and for any prime $\fP \in S|_M$ above a prime $\fp \in S|_L$.
\end{proof}
There is however a slight caveat in working with units modulo torsion, because if $F/k$ is a Galois
extension with Galois group $G$, then the fixed submodule of $U_S(F)$ under a subgroup $H$
of $G$ need not be canonically isomorphic to $U_S(F^H)$. We will need to understand
exactly when it is and what the difference is whenever it is not. Write $\mu(M)$ for the group of roots
of unity of a number field $M$. Then, from the short exact sequence
\[1\rightarrow \mu(F)\rightarrow \Sunits{F}\rightarrow U_S(F) \rightarrow 1\]
we get the usual long exact sequence of group cohomology
\[1\rightarrow \mu(F)^H\rightarrow (\Sunits{F})^H\rightarrow U_S(F)^H\rightarrow\sH^1(H,\mu(F))\rightarrow
\sH^1(H,\Sunits{F})\]
for any subgroup $H$ of $G$. We see immediately that $(\Sunits{F})^H/\mu(F)^H=\Sunits{F^H}/\mu(F^H)$ is
canonically isomorphic to $U_S(F)^H$ if and only if $\ker(\sH^1(H,\mu(F))\rightarrow \sH^1(H,\Sunits{F}))=0$.
We have that $f\in \sH^1(H,\mu(F))$ is in this kernel if and only if there is an $S$-unit $u\in \Sunits{F}$
such that $f(h) = h(u)/u\in \mu(F)\;\forall h\in H$. If $f$ is not a co-boundary itself, then $u\notin \mu(F)$
and $u$ can without loss of generality be taken to be non-torsion. We deduce that $F$ must contain a
root of a non-torsion $S$-unit of an intermediate extension of $F/F^H$. Conversely, if it does, then defining
$f$ as above gives a non-trivial element of the kernel. In summary, we record
\begin{lem}\label{lem:index}
With $F/k$ and $S$ as above, we have, for any subgroup $H$ of the Galois group of $F/k$, that
$U_S(F)^H\cong U_S(F^H)$ if no intermediate extension of $F/F^H$ is obtained by
adjoining a root of a non-torsion $S$-unit. In general,
\[[U_S(F)^H:U_S(F^H)] = \#\ker(\sH^1(H,\mu(F))\rightarrow \sH^1(H,\Sunits{F})).\]
\end{lem}
We are now ready to prove the main result of this section, which links regulator quotients to Dokchitser
constants.
\begin{prop}\label{prop:newReg}
Let $F/k$ be a finite Galois extension with Galois group $G$, let $S$ be a finite $G$-stable set of
places of $F$ including the Archimedean ones.
Write $\lambda(H)= \#\ker(\sH^1(H,\mu(F))\rightarrow \sH^1(H,\Sunits{F}))$.
If $\Theta=\sum_i H_i - \sum_j H_j'$ is a $G$-relation, then we have
\begin{eqnarray*}
\C_\Theta(U_S(F)) = \frac{\C_\Theta(\triv)}{\prod\limits_{\fp\in S|_k}\C_\Theta(\bZ[G/D_\fp])}\cdot
\frac{\prod_i \left(R_S(F^{H_i})/\lambda(H_i)\right)^2}{\prod_j \left(R_S(F^{H_j'})/\lambda(H_j')\right)^2},
\end{eqnarray*}
where, for each $\fp\in S|_k$, $D_\fp$ is a decomposition subgroup of $G$ at $\fp$.
\end{prop}
\begin{rmrk}
In particular, if the decomposition groups at all primes in $S$ are cyclic, which is for example the
case if $S$ is the set of Archimedean places, then the associated Dokchitser constants are trivial
by \cite[Lemma 2.46]{TVD-2} and the formula just reads
\[\C_\Theta(U_S(F)) = \C_\Theta(\triv)\cdot
\frac{\prod_i \left(R_S(F^{H_i})/\lambda(H_i)\right)^2}{\prod_j \left(R_S(F^{H_j'})/\lambda(H_j')\right)^2}.\]
\end{rmrk}
\begin{proof}
For any $H\leq G$ we have
\begin{eqnarray*}
\lambda(H)^2\cdot\det\left(\frac{1}{|H|}\left<,\right>_F\big|U_S(F)^H\right) &
\stackrel{\text{Lemma }\ref{lem:index}}{=} & \det\left(\frac{1}{|H|}\left<,\right>_F\big|U_S(F^H)\right)\\
& \stackrel{\text{Lemma }\ref{lem:baseChange}}{=} & \det\left(\left<,\right>_{F^H}\big|U_S(F^H)\right)\\
& \stackrel{\text{Lemma }\ref{lem:detPairing}}{=} & 
\frac{\sum_{\fP\in S|_{F^H}}f_\fP e_\fP}{\prod_{\fP\in S|_{F^H}}f_\fP e_\fP}\cdot R_S(F^H)^2.
\end{eqnarray*}
Since $S$ contains precisely all the places above the places in $S|_k$, we have, for each $\fp\in S|_k$,
\[\sum_{\substack{\fP\in S|_{F^H},\\ \fP|\fp}} f_\fP e_\fP =
\sum_{\substack{\fP\in S|_{F^H},\\ \fP|\fp}} f_\fp e_\fp f_{\fP/\fp}e_{\fP/\fp} = f_\fp e_\fp [F^H:k]\]
and thus, for any
$H\leq G$ we have
\[\sum_{\fP\in S|_{F^H}} f_\fP e_\fP = [F^H:k]\sum_{\fp\in S|_k}f_\fp e_\fp = 
\frac{|G|}{|H|}\sum_{\fp\in S|_k}f_\fp e_\fp.\]
Observe that by \cite[Example 2.30]{TVD-2} the term $\sum_{\fp\in S|_k}e_\fp f_\fp$, being a constant,
vanishes in a relation.

Also, for each $\fp\in S|_k$ we have
\[\prod_{\substack{\fP\in S|_{F^H},\\ \fP|\fp}}f_\fP e_\fP =  (f_\fp e_\fp)^{\#\{\text{primes of }F^H\text{ above }\fp\}}
\cdot\prod_{\substack{\fP\in S|_{F^H},\\ \fP|\fp}}(f_{\fP/\fp} e_{\fP/\fp}).\]
Now, by \cite[Example 2.37]{TVD-2} the first of the two factors vanishes in
a relation. Also, by \cite[Corollary 2.44]{TVD-2}, the second
factor may be replaced by $\C_{\Theta}(\bZ[G/D])$ in a relation, where $D$ is a decomposition subgroup
of $G$ at $\fp$. Combining everything we have said, we obtain
\begin{eqnarray*}
\C_\Theta(U_S(F)) & = &\frac{\prod_i \det\left(\frac{1}{|H_i|}\left<,\right>_F\big|U_S(F)^{H_i}\right)}{
\prod_j \det\left(\frac{1}{|H_j'|}\left<,\right>_F\big|U_S(F)^{H_j'}\right)}\\
& = & \frac{\C_\Theta(\triv)}{\prod\limits_{\fp\in S|_k}\C_\Theta(\bZ[G/D_\fp])}\cdot
\frac{\prod_i \left(R_S(F^{H_i})/\lambda(H_i)\right)^2}{\prod_j \left(R_S(F^{H_j'})/\lambda(H_j')\right)^2},
\end{eqnarray*}
as claimed.
\end{proof}
Now that we know how to turn quotients of regulators of number fields into Dokchitser constants,
which are purely representation theoretic invariants,
we will establish some properties of Dokchitser constants in the next two sections.
\section{An alternative description of Dokchitser constants}\label{sec:altDefn}
The definition of Dokchitser constants that we have given above is somewhat unsatisfactory, since
it involves making an arbitrary choice (that of a pairing) on which the result does not depend. It
would be nice to have a definition that avoids any arbitrary choices. As a first step in the 
investigation of the properties of Dokchitser constants, we will provide an alternative definition, 
which depends on fixing more specific information about the
relation (on which the result again does not depend) but not on any choices connected with the
representation.

Let $\Theta = \sum_i H_i - \sum_j H_j'$ be an $AG$-relation. Define the $G$-sets
$S_1=\bigsqcup_i G/H_i$ and $S_2=\bigsqcup_j G/H_j'$. Then to say that $\bQ[S_1] \cong \bQ[S_2]$ is
equivalent to saying that there exists an embedding of $\bZ G$-lattices
$$\phi:\bZ[S_1]\hookrightarrow\bZ[S_2]$$
with finite cokernel. Also, to say that $\bZ_{(l)}[S_1] \cong \bZ_{(l)}[S_2]$ is equivalent to
saying that there is such a $\phi$ with finite cokernel of order coprime to $l$.

With these remarks in mind, let $\cR$ be a PID containing $\bZ$, let $\Gamma$ be an $\cR G$-lattice
and fix an injection of $\cR G$-modules
$$\phi:\cR[S_1]\hookrightarrow \cR[S_2].$$

By dualising this, we obtain
$$\phi^{\tr}:\Hom(\cR[S_2],\cR)\hookrightarrow \Hom(\cR[S_1],\cR).$$
Recall that if $\Gamma$ is any $\cR G$-lattice
and if $g\in G$ is represented by the matrix $M$ with respect to some $\cR$-basis,
then in the action on $\Hom(\Gamma,\cR)$, $g$ is represented by $(M^{-1})^\tr$ with
respect to the dual basis. But this is equal to $M$ if $M$ is a permutation matrix
and so permutation modules are canonically self-dual. In summary, we have a map
$$\phi^\tr:\cR[S_2]\hookrightarrow \cR[S_1].$$
Applying the functor $\Hom(\_\;,\Gamma)$ yields
$$\phi^*:\Hom_\cR(\cR[S_2],\Gamma)\rightarrow \Hom_\cR(\cR[S_1],\Gamma)$$
and
$$(\phi^\tr)^*:\Hom_\cR(\cR[S_1],\Gamma)\rightarrow \Hom_\cR(\cR[S_2],\Gamma).$$
Upon restricting to the $G$-invariant subspaces we obtain maps $\phi_G^*$ and $(\phi^\tr)_G^*$
between the corresponding spaces of $G$-homomorphisms (to avoid index overload,
we are abusing notation slightly by not including $\Gamma$ in the notation of these maps). Since
$\cR$ is a PID, the spaces of $G$-homomorphisms are $\cR$-free. Also, since $\phi\otimes \cK$ and
$\phi^\tr\otimes \cK$ are both isomorphisms, so are $\phi^*_G\otimes \cK$ and
$(\phi^\tr)^*_G\otimes \cK$. Thus both $\phi_G^*$ and $(\phi^\tr)^*_G$ have non-zero determinants.
\begin{defn}
Define
$$\C_\Theta'(\Gamma)=\frac{\det(\phi^\tr)_G^*}{\det\phi_G^*}\in \cK^\times/(\cR^\times)^2$$
with both determinants computed with respect to the same bases on $\Hom_{\cR[G]}(\cR[S_i],\Gamma)$, $i=1,2$. If we change the basis on $\Hom_{\cR[G]}(\cR[S_1],\Gamma)$, say, then the quotient
changes by the square of the determinant of change of basis, so it really is a well-defined element
of $\cK^\times/(\cR^\times)^2$.
\end{defn}
The injection $\phi$ is not present in the notation and indeed it will turn out that the residue
class of $\C_\Theta'(\Gamma)$ modulo $(\cR^\times)^2$ only depends on the relation and the module
but not on any other choices. The main result of this section is
\begin{thm}\label{thm:altDefn}
Let $G$ be a finite group, $\cR$ a principal ideal domain containing $\bZ$ with field of
fractions $\cK$, $\Theta=\sum_i H_i - \sum_j H_j'$ an
$AG$-relation, where $A$ is either $\bQ$ or $\bZ_{(p)}$, and $\Gamma$ an $\cR G$-lattice. Fix an
injection $\phi:\cR[S_1]\hookrightarrow \cR[S_2]$ and obtain $\phi_G^*$ and $(\phi^\tr)_G^*$ as above.
Fix a $G$-invariant non-degenerate bilinear pairing $\left<,\right>$ on $\Gamma$
(which exists because $\Gamma\otimes \cK$ is self-dual by convention). Then
$$\frac{\det(\phi^\tr)_G^*}{\det\phi_G^*} \equiv
\frac{\prod_i \det\left(\frac{1}{|H_i|}\left<,\right>\big|\Gamma^{H_i}\right)}
{\prod_j \det\left(\frac{1}{|H'_j|}\left<,\right>\big|\Gamma^{H'_j}\right)} \mod (\cR^\times)^2.$$
\end{thm}
\begin{proof}
Define a pairing $(,)_1$ on $\Hom_{\cR}(\cR[S_1],\Gamma)$ by
$$(f_1,f_2)_1 = \frac{1}{|G|}\sum_{s\in S_1}\left<f_1(s),f_2(s)\right>$$
and define an analogous pairing $(,)_2$ on $\Hom_{\cR}(\cR[S_2],\Gamma)$. It is immediate that these
pairings, when restricted to the spaces of $G$-homomorphisms, are $G$-invariant. We first claim that
$(\phi^\tr)^*$ is the adjoint of $\phi^*$ with respect to these pairings. Indeed, it suffices to
check this for any particular choice of bases on $\Hom_\cR(\cR[S_1],\Gamma)$ and on $\Hom_\cR(\cR[S_2],\Gamma)$.
Let $S_1=\left\{s_1,\ldots,s_n\right\}$ and choose a basis $\gamma_j$, $j=1,\ldots,r$ of $\Gamma$.
Define $f_{i,j}\in \Hom_\cR(\cR[S_1],\Gamma)$ by $f_{i,j}(s_i) = \gamma_j$, $f_{i,j}(s)=0\;\forall s\neq s_i$.
Then $f_{i,j},\;i=1,\ldots,n,\;j=1,\ldots,r$ is a basis of $\Hom_\cR(\cR[S_1],\Gamma)$. Fix the analogous
basis $f_{i,j}'$ for $\Hom_\cR(\cR[S_2],\Gamma)$ where $S_2=\left\{s_1',\ldots,s_n'\right\}$.
We compute
\begin{eqnarray}\label{eq:adjoint}
|G|\cdot(f_{i,j},\phi^*f_{r,t}')_1 & = & \sum_{s\in S_1} \left<f_{i,j}(s),f_{r,t}'(\phi(s))\right>
= \left<\gamma_j,f_{r,t}'(\phi(s_i))\right>\nonumber\\
& = & \left<\gamma_j,\phi_{i,r}\gamma_t\right>
= \left<\phi_{i,r}\gamma_j,\gamma_t\right>
= \left<f_{i,j}(\phi^\tr s_r'),\gamma_t\right>\nonumber\\
& = & \sum_{s'\in S_2} \left<f_{i,j}(\phi^\tr s'),f_{r,t}'(s')\right>
=|G|\cdot((\phi^\tr)^*f_{i,j},f_{r,t}')_2
\end{eqnarray}
as required. Next, for a subgroup $H$ of $G$ we can identify $\Hom_G(G/H,\Gamma)$ with $\Gamma^H$
via $f\mapsto f(1)$. We claim that under this identification, we have
\begin{eqnarray}\label{eq:newPairing}
\det\left((,)_1|\Hom_{\cR[G]}(\cR[S_1],\Gamma)\right) \equiv
\prod_i \det\left(\frac{1}{|H_i|}\left<,\right>\big|\Gamma^{H_i}\right)\mod (\cR^\times)^2
\end{eqnarray}
and similarly for $S_2$. Indeed, if for subgroups $H_i\neq H_k$, we have that $\cR[G/H_i]$ and $\cR[G/H_k]$
are summands of $\cR[S_1]$, then an element of $\Hom_{\cR[G]}(\cR[S_1],\Gamma)$ which
is trivial outside of $G/H_i$ is orthogonal to an element which is trivial outside of $G/H_k$. So it
suffices to prove the claim for $S_1=G/H$. We compute
\begin{eqnarray*}
(f_1,f_2)_1 & = & \frac{1}{|G|}\sum_{s\in G/H} \left<f_1(s),f_2(s)\right>
= \frac{1}{|G|}\sum_{s\in G/H} \left<s\cdot f_1(1),s\cdot f_2(1)\right>\\
& = & \frac{1}{|G|}\sum_{s\in G/H} \left<f_1(1),f_2(1)\right>
= \frac{1}{|H|}\left<f_1(1),f_2(1)\right>,
\end{eqnarray*}
which immediately implies the claim. Now, fix bases $v_1,\ldots,v_m$ and $v_1',\ldots,v_m'$ on
$\Hom_{\cR[G]}(\cR[S_1],\Gamma)$ and $\Hom_{\cR[G]}(\cR[S_2],\Gamma)$, respectively. Then
\begin{eqnarray*}\frac{\prod_i \det\left(\frac{1}{|H_i|}\left<,\right>\big|\Gamma^{H_i}\right)}
{\prod_j \det\left(\frac{1}{|H'_j|}\left<,\right>\big|\Gamma^{H'_j}\right)} &
\stackrel{\text{by } (\ref{eq:newPairing})}{\equiv} &
\frac{\det\left((v_i,v_j)_1\right)}
{\det\left((v_k',v_l')_2\right)}\\
& \equiv & \frac{\det\left((v_i,\phi_G^*v_l')_1\right)/\det(\phi_G^*)}
{\det\left(((\phi^\tr)_G^*v_i,v_l')_2\right)/\det((\phi^\tr)_G^*)}\\
& \stackrel{\text{by } (\ref{eq:adjoint})}{\equiv} & \det((\phi^\tr)_G^*)/\det(\phi_G^*)\mod (\cR^\times)^2,
\end{eqnarray*}
which concludes the proof.
\end{proof}
\begin{cor}
The value of $\C_\Theta'(\Gamma)$ is independent of the choice of $\phi$ and we have
$\C_\Theta'(\Gamma) = \C_\Theta(\Gamma)$ for all $G$-relations $\Theta$ and all $\cR G$-lattices
$\Gamma$.
\end{cor}
\begin{rmrk}
It is not difficult to prove the independence of $\phi$ directly, using explicit calculations with
the bases $f_{i,j}:s_l \mapsto \delta_{i,l}\cdot\gamma_j$ and $f_{i,j}':s_l' \mapsto \delta_{i,l}\cdot\gamma_j$ from the above proof. This also gives an alternative
proof of the independence of $\C_\Theta(\Gamma)$ of the pairing.
\end{rmrk}
\begin{rmrk}
It is instructive to compare our alternative definition of Dokchitser constants in conjunction with
Proposition \ref{prop:newReg} with the formula for
class number quotients derived by de Smit in \cite[Theorem 2.2]{Smi-01}. In his formula, the torsion
subgroup of the units is more directly incorporated into the whole expression. However, arbitrary
$\bZ G$-modules of a given group $G$ are more difficult to classify than $\bZ$-free modules and we will need to use
the classification from \cite{Lee-64} for $G=D_{2p}$ in the next section. That is the main reason why
we pass to the quotient modulo torsion first and then recover the torsion separately in the shape
of $\lambda(H)$. Another reason to work with Dokchitser constants is that a lot of the computations
of Dokchitser constants in the next section will be easier using a pairing rather than an embedding $\phi$.
\end{rmrk}
An immediate consequence of the alternative definition is the following:
\begin{lem}\label{lem:zpRelations}
Let $G$ be a finite group and $\Theta$ a $\bZ_{(l)} G$-relation. Then $\ord_l(\C_\Theta(\Gamma))=0$
for all $\bZ G$-lattices~$\Gamma$.
\end{lem}
\begin{proof}
As remarked above, we can find an injection of $G$-modules $\phi:\bZ[S_1]\hookrightarrow \bZ[S_2]$ with
co-kernel of order coprime to $l$. For any injection of free abelian groups with finite co-kernel,
the order of its co-kernel is equal to the absolute value of its determinant (with respect to any bases).
Now, applying $\Hom_G(\_\;,\Gamma)$ to the tautological short exact sequence
\[0\rightarrow \bZ[S_1]\stackrel{\phi}{\rightarrow} \bZ[S_2]\rightarrow \coker(\phi)\rightarrow 0\]
gives the long exact sequence
\[ 0\rightarrow \Hom_G(\bZ[S_2],\Gamma)\stackrel{\phi_G^*}{\rightarrow} \Hom_G(\bZ[S_1],\Gamma)
\rightarrow \Ext_{\bZ G}^1(\coker(\phi),\Gamma).\]
Since $\Ext_{\bZ G}^1(\coker(\phi),\Gamma)$ has no $l$-torsion, neither does $\coker(\phi^*_G)$.
The same goes for $(\phi^\tr)^*_G$ and the proof is complete.
\end{proof}
\begin{defn}
Let $l$ be a prime number. A finite group is called $l$-\emph{hypo-elementary} if it has a normal Sylow $l$-subgroup with cyclic
quotient. Equivalently, an $l$-hypo-elementary group is a semi-direct
product of an $l$-group acted on by a cyclic group of order coprime to $l$.
\end{defn}
\begin{thm}[Conlon's Induction Theorem]\label{thm:Conlon}
Given any finite group $H$ and any
commutative ring $\tilde{R}$ in which every prime divisor of $|H|$ except possibly $l$ is invertible,
there exist integers $\alpha_{H'}$ such that some integer multiple of the trivial representation of
$H$ over $\tilde{R}$ is equal to $\sum_{H'} \alpha_{H'} \tilde{R}[H/H']$ in the representation ring
over $\tilde{R}$, where the sum is taken over $l$-hypo-elementary subgroups of $H$.
\end{thm}
A proof can be found e.g. in \cite{CR}, (80.60). We will use this result to considerably strengthen
Proposition \ref{prop:notDivble}:
\begin{prop}\label{prop:normalSgpCycQuo}
Let $G$ be a finite group, let $N$ be a normal subgroup such that the quotient
group $C=G/N$ is cyclic. Let $l$ be a prime not dividing the order of $N$ and let $\cR=\bZ$ or
$\bZ_l$. Then
$$\ord_l(\C_\Theta(\Gamma))=0$$
for all $\cR G$-lattices $\Gamma$ and all $G$-relations $\Theta$.
\end{prop}
\begin{proof}
By Lemma \ref{lem:zpRelations}, it suffices to show that every $\bQ G$-relation is in fact a
$\bZ_{(l)} G$-relation. For that, it is enough to show that the rank of the sublattice of
$\bZ_{(l)} G$-relations is equal to the rank of the lattice of $\bQ G$-relations, since the former
is saturated in the latter.\footnote{The lattice of $\bZ_{(l)} G$-relations is the kernel of the
natural map from the Burnside ring to the representation ring over $\bZ_{(l)}$. This map is clearly
linear and kernels of linear maps from abelian groups are always saturated.} By Theorem \ref{thm:Conlon}, the rank of the lattice of
$\bZ_{(l)} G$-relations is at least equal to the number of conjugacy classes of
non-$l$-hypo-elementary subgroups. Explicitly,
for each subgroup $H$ of $G$ which is not $l$-hypo-elementary, we get a $\bZ_{(l)} G$-relation
$\alpha_H H - \sum_{H'} \alpha_{H'}H'$, the sum taken over $l$-hypo-elementary subgroups of $H$.
All relations obtained in this way are clearly linearly independent, since each one contains a
unique 'maximal' subgroup that has the property that all other subgroups featuring in the relation
are contained in this one.
Since the rank of the lattice of $\bQ G$-relations is equal to the number of conjugacy classes of
non-cyclic subgroups of $G$, it is enough to show that any $l$-hypo-elementary subgroup of $G$ must
by cyclic.

So take $H=P\rtimes Z\leq G$ where $P$ is an $l$-group and $Z$ is cyclic of order
coprime to $l$.
Since $l$ does not divide $|N|$ we have that
\[P \cong P/P\cap N \cong PN/N \leq G/N\]
is cyclic. Further, since $H/P$ is abelian, the commutator subgroup $H'$ of $H$ must lie in $P$, so
it is an $l$-group. But also, $H'\leq G'\leq N$, since $G/N$ is abelian. Therefore $H'=\{1\}$, since
$l$ does not divide $|N|$. It follows that $H$ is abelian, $H = P \times C$ and so cyclic.
\end{proof}
\section{Dokchitser constants in dihedral groups}\label{sec:dihGroups}
We will now compute the Dokchitser constants of all $\bZ G$-lattices when
$G=D_{2p}$ is the dihedral group with $2p$ elements for $p$ an odd prime and $\Theta$ is the relation
from Example \ref{ex:d2p}. By Proposition \ref{prop:multiplicativity}, we only need to compute them
for indecomposable representations. Nonetheless, the fact that this can be done at all is a piece
of good fortune. We will begin by recalling the classification of indecomposable integral
representations of $D_{2p}$ from \cite{Lee-64}.

Fix a primitive $p$-th root of unity $\zeta_p$ in a fixed algebraic closure $\bar{\bQ}$ of $\bQ$.
Let $\bQ(\z)^+$ be the maximal real subfield of the $p$-th cyclotomic field and let ${\eusm O}^+$
be its ring of integers. Let $\{U_i\}$ be a full set of representatives of the ideal class group of
$\bQ(\zeta_p)^+$ and take $U_1 = U = {\eusm O}^+$ to represent the principal ideals. Write
$G = \left<\sigma,\varpi:\;\sigma^2 = \varpi^p = (\sigma\varpi)^2 = 1\right>$. Let ${\eusm O}$ be the
ring of integers of
$\bQ(\z)$. Write $A_i$ for the $\bZ G$-module $U_i{\eusm O}$ on which $\sigma$ acts as complex
conjugation and $\varpi$ as multiplication by $\z$. Let $A_i'$ be the module
$(\bar{\z}-\z)U_i{\eusm O}$
with the same $G$-action. Set $A=A_1$, $A'=A_1'$.

Finally write $\triv$ for the 1-dimensional trivial $\bZ G$-module, $\epsilon$ for the 1-dimensional
module sending $\sigma$ to -1 and $\varpi$ to 1 and $\Delta$ for the 2-dimensional module
$\bZ[G/C_p]$ which is an extension of $\triv$ by $\epsilon$. The following is a complete list of
non-isomorphic indecomposable $\bZ G$-lattices (see \cite{Lee-64}):
\begin{itemize}
	\item $\triv$;
	\item $\epsilon$;
	\item $\Delta$;
	\item for each $i$, $A_i$;
	\item for each $i$, $A_i'$;
	\item for each $i$, a non-trivial extension of $\triv$ by $A_i'$, denoted by $(A_i',\triv)$;
	\item for each $i$, a non-trivial extension of $\epsilon$ by $A_i$, denoted by $(A_i,\epsilon)$;
	\item for each $i$, a non-trivial extension of $\Delta$ by $A_i$, denoted by $(A_i,\Delta)$;
	\item for each $i$, a non-trivial extension of $\Delta$ by $A_i'$, denoted by $(A_i',\Delta)$;
	\item for each $i$, a non-trivial extension of $\Delta$ by $A_i\oplus A_i'$,
	denoted by $(A_i\oplus A_i',\Delta)$;
\end{itemize}
It is a trivial check that $\C_\Theta(\triv) = 1/p$, $\C_\Theta(\epsilon) = p$,
$\C_\Theta(\Delta) = 1$.
\begin{lem}\label{lem:regConsts1}
The Dokchitser constants of $A$ and of $A'$ are $p$ and $1/p$, respectively.
\end{lem}
\begin{proof}
The matrices of $\sigma$, $\varpi$ acting on $A'$ on the left with respect to the basis
$\bar{\z}-\z, (\bar{\z}-\z)\z,(\bar{\z}-\z)\z^2,\ldots,(\bar{\z}-\z)\z^{p-2}$ are 
$$\left[  
\begin{array}{c c c c c c c}
	-1 & 1 & 0 & 0 & 0 & \cdots & 0 \\
	0 & 1 & 0 & 0 & 0 & \cdots & 0 \\
	0 & 1 & 0 & 0 & \cdots & 0 & -1 \\
	0 & 1 & 0 & \cdots & 0 & -1 & 0 \\
	\vdots & \vdots & \vdots & \adots & \adots & \adots  & \vdots \\
	0 & 1 & 0 & -1 & 0 & \cdots & 0 \\
	0 & 1 & -1 & 0 & \cdots & 0 & 0 \\
\end{array}
\right]
\text{ and }
\left[ \begin{array}{c c c c c c}
	0 & 0 & 0 & \ldots & 0 & -1 \\
	1 & 0 & 0 & \ldots & 0 & -1 \\
	0 & 1 & 0 & \ldots & 0 & -1 \\
	\vdots & \vdots & \ddots & \ddots & \vdots & \vdots \\
	0 & \ldots & 0 & 1 & 0 & -1 \\
	0 & \ldots & 0 & 0 & 1 & -1 \\
\end{array} \right],
$$
respectively. It is immediately seen that the same matrices represent the $G$-action by multiplication
on the submodule
$$\left<\varpi^{\frac{p-1}{2}}-\varpi^{\frac{p+1}{2}},\varpi^{\frac{p+1}{2}}-\varpi^{\frac{p+3}{2}},
\ldots,\varpi^{p-1}-1,1-\varpi,
\ldots,\varpi^{\frac{p-5}{2}}-\varpi^{\frac{p-3}{2}}\right>_\bZ$$
of $\bZ[G/C_2]$ with respect to the indicated basis. But this is just the submodule
$$\left<1-\varpi^i: i\in \{1,\ldots,p-1\}\right>_\bZ$$
of the permutation lattice $\bZ[G/C_2]$. We can choose the standard pairing on the latter which makes
the different coset representatives an orthonormal $\bZ$-basis. It is easy to see that the fixed 
sublattices under 1 and under $\left<\sigma\right> = C_2$ are
$$\left<1-\varpi^i: i=1,\ldots,p-1\right>_\bZ \text{ and } \left<2-\varpi^i-\varpi^{p-i}: i = 1,
\ldots,\frac{p-1}{2}\right>_\bZ,$$
respectively. The subgroup $C_p$ only fixes the trivial lattice. The matrices of the pairing on these 
modules with respect to the bases indicated are then
$$\left[ \begin{array}{c c c c c}
	2 & 1 & 1 & \cdots & 1 \\
	1 & 2 & 1 & \cdots & 1 \\
	\vdots & \ddots & \ddots & \ddots & \vdots \\
	1 & \cdots & 1 & 2 & 1 \\
	1 & 1 & \cdots & 1 & 2 \\
\end{array}\right]
\text{ and }
\left[ \begin{array}{c c c c c}
	6 & 4 & 4 & \cdots & 4 \\
	4 & 6 & 4 & \cdots & 4 \\
	\vdots & \ddots & \ddots & \ddots & \vdots \\
	4 & \cdots & 4 & 6 & 4 \\
	4 & 4 & \cdots & 4 & 6 \\
\end{array}\right]$$
of sizes $p-1$ and $\frac{p-1}{2}$ with determinants $p$ and $2^\frac{p-1}{2}p$, respectively, as
can be checked by 
elementary row operations. So, taking into account the normalisation by the sizes of the subgroups,
we get
$$\frac{\text{det}\left(\frac{1}{|1|}\left<,\right>|A'^{1}\right)\text{det}\left(\frac{1}{|G|}\left<,
\right>|A'^{G}\right)^2}
{\text{det}\left(\frac{1}{|C_2|}\left<,\right>|A'^{C_2}\right)^2\text{det}\left(\frac{1}{|C_p|}\left<,
\right>|A'^{C_p}\right)}
= \frac{p}{p^2} = 1/p$$
as claimed.

Now consider the $\bZ G$-module $\bZ[G/C_2] \otimes_\bZ \epsilon$ with diagonal $G$-action. It is now
clear from above that $A$ is isomorphic to the submodule of $\bZ[G/C_2] \otimes_\bZ \epsilon$ given by
$\left<1-\varpi^i: i = 1,\ldots,p-1\right>$. The fixed submodules under $1$ and under $C_2$ are
$$\left<1-\varpi^i: i = 1,\ldots,p-1\right>\text{ and }\left<\varpi^i-\varpi^{p-i}: i =
1\ldots \frac{p-1}{2}\right>,$$
respectively, and an entirely similar calculation using the same natural pairing as above shows that
$\C_\Theta(A) = p$.
\end{proof}
\begin{lem}\label{lem:regConsts2}
We have $(A',\triv)\cong \bZ\left[G/C_2\right]$ and ${\C}_\Theta((A',\triv)) = 1$.
\end{lem}
\begin{proof} Take the $\bZ$-basis $1,\varpi,\ldots,\varpi^{p-1}$ for $\bZ\left[G/C_2\right]$.
Then there
is the submodule 
$\left<\sum_{i=0}^{p-1}\varpi^i\right>$
isomorphic to $\triv$ and the submodule
$\left<1-\varpi^i: i\in \{1,\ldots,p-1\}\right>$
isomorphic to $A'$ but their sum is the submodule $\left\{\sum_i \alpha_i\varpi^i: \sum \alpha_i
\equiv 0 (\text{ mod }p)\right\}$ which is an index $p$ sublattice. In fact $\Gamma=\bZ[G/C_2]$ is
indecomposable,
for if it were decomposable, it would have to decompose as
$\widehat{\Gamma_1}\cap\Gamma\oplus\widehat{\Gamma_2}\cap\Gamma$, where
$\Gamma\otimes\bQ=\widehat{\Gamma_1}\oplus\widehat{\Gamma_2}$ is the decomposition into irreducible
rational representations. But these intersections are easily seen to be the sublattices just
exhibited.
Thus $\bZ[G/C_2]$ must be a non-trivial extension of $\triv$ by $A'$ and the first claim follows
from the classification of integral representations. The Dokchitser constant of $(A',\triv)$ is 
therefore trivial by \cite[Lemma 2.46]{TVD-2}.
\end{proof}

\begin{lem}\label{lem:regConsts3}
The Dokchitser constants of the remaining lattices in the above list for $i=1$ are as follows:
\begin{itemize}
	\item ${\C}_\Theta((A,\epsilon)) = 1$;
	\item ${\C}_\Theta((A,\Delta)) = 1/p$;
	\item ${\C}_\Theta((A',\Delta)) = p$;
	\item ${\C}_\Theta((A\oplus A',\Delta)) = 1$;
\end{itemize}
\end{lem}
\begin{proof}
It is noted in \cite[\S 4]{Lee-64} that $(A\oplus A',\Delta)\cong \bZ[G/1]$ and so ${\C}_\Theta((A\oplus A',\Delta)) = 1$ by \cite[Lemma 2.46]{TVD-2}.

For the other three lattices, since we only need to determine the $p$-parts it suffices to work up to squares of elements with trivial $p$-valuation so we will work over $\bZ_p$ rather than over $\bZ$. Write $\widetilde{(A,\epsilon)} = (A,\epsilon)\otimes_\bZ\bZ_p$ and similarly for the other lattices. Since $\triv\oplus\epsilon$ is an index 2 sublattice of $\Delta$, over $\bZ_p$ we have $\tilde{\triv}\oplus\tilde{\epsilon}\cong \tilde{\Delta}$. Now, $(A,\epsilon)\otimes\epsilon\cong (A',\triv)$ and so
\begin{eqnarray*}
\widetilde{(A,\epsilon)}\oplus\widetilde{(A',\triv)} & \stackrel{\ref{lem:regConsts2}}{=} & \bZ_p[G/C_2]\otimes(\tilde{\triv}\oplus\tilde{\epsilon})\\
 & \cong & \bZ_p[G/C_2]\otimes\tilde{\Delta} \\
 & \cong & \bZ_p[G/C_2]\otimes\bZ_p[G/C_p] \\
 & \cong & \bZ_p[G/1]
\end{eqnarray*}
which has trivial Dokchitser constant by \cite[Lemma 2.46]{TVD-2}. By multiplicativity of Dokchitser constants and by Lemma \ref{lem:regConsts2}, ${\C}_\Theta(\widetilde{(A,\epsilon)}) = 1$. Similarly,
$\widetilde{(A,\Delta)}\cong (\tilde{A},\tilde{\triv}\oplus\tilde{\epsilon})$, and since
$\Ext_{\bZ G}^1(\triv,A)=0$ (\cite[Lemma 2.1]{Lee-64}), it is easy to see that
$$(\tilde{A},\tilde{\triv}\oplus\tilde{\epsilon}) \cong \tilde{\triv}\oplus \widetilde{(A,\epsilon)}.$$
By multiplicativity of Dokchitser constants, we deduce that
\[{\C}_\Theta(\widetilde{(A,\Delta)}) = 1/p \in \bQ_p^\times/\left(\bZ_p^\times\right)^2.\]
Also $\Ext_{\bZ G}^1(\epsilon,A')=0$ and
\[(\tilde{A'},\tilde{\triv}\oplus\tilde{\epsilon}) \cong \tilde{\epsilon}\oplus \widetilde{(A',\triv)},\]
whence
\[{\C}_\Theta(\widetilde{(A,\Delta)}) = p \in \bQ_p^\times/\left(\bZ_p^\times\right)^2.\]
\end{proof}
\begin{thm}\label{thm:allRegConsts}
The Dokchitser constants of all the indecomposable integral representations of $D_{2p}$ for $p$ an
odd prime are as follows:
\begin{center}
\begin{tabular}{| r | l | l | l | l | l | l | l | l | l | l |}
\hline
$\Gamma$: & $\triv$ & $\epsilon$ & $\Delta$ & $A_i$ & $A_i'$ & $(A_i',\triv)$ & $(A_i,\epsilon)$ & $(A_i,\Delta)$ & $(A_i',\Delta)$ & $(A_i\oplus A_i',\Delta)$\\
\hline
${\C}_\Theta(\Gamma)$: & $1/p$ & $p$ & $1$ & $p\;\forall i $ & $1/p\;\forall i $ &
$1\;\forall i $ & $1\;\forall i $ & $1/p\;\forall i $ & $p\;\forall i $ & $1\;\forall i $\\
\hline
\end{tabular}
\end{center}

\end{thm}
\begin{proof}
For $i=1$ this is Lemma \ref{lem:regConsts1}, Lemma \ref{lem:regConsts2} and Lemma
\ref{lem:regConsts3}. We will show that $A_i$ is isomorphic to $A$ over $\bZ_{(2)}$ and over
$\bZ_{(p)}$ for all $i$ and $A_i'$ is isomorphic to $A'$ over $\bZ_{(2)}$ and over $\bZ_{(p)}$ for
all $i$ (strictly speaking, the isomorphism over
$\bZ_{(p)}$ would be enough for this theorem by Proposition \ref{prop:normalSgpCycQuo}).
Recall that $A_i$, $A_i'$ are given by $(\bar{\z}-\z)^jU_i{\eusm O}$ for $j=0,1$, respectively,
where $U_i$ runs through representatives of the ideal class group of $\bQ(\z)^+$. Take each $U_i$
to be of norm coprime to $2p$. Then $A_i$ is a sublattice of $A=A_1$ of index coprime to $2p$
and the two are therefore isomorphic over $\bZ_2$ and over $\bZ_p$. Thus they have the same
Dokchitser constants. Similarly, $A_i'$ all have the same Dokchitser constants as $A'=A_1'$.
\end{proof}
The proof of Theorem \ref{thm:allRegConsts} exhibits an important feature of Dokchitser constants,
which we will now summarise.
\begin{defn}
Given a finite group $G$ and a principal ideal domain $\cR$, two $\cR G$-lattices $M$ and $N$ are
said to lie in the same \emph{genus} if $M\otimes \cR_{\mathfrak p}\cong N\otimes \cR_{\mathfrak p}$
as $\cR_{\mathfrak p}G$-modules for all completions $\cR_{\mathfrak p}$ at prime ideals
${\mathfrak p}$ of $\cR$. This is clearly an equivalence relation.
\end{defn}
The two main conceptual steps in the proof of Theorem \ref{thm:allRegConsts} can be summarised as:
\begin{thm}
The Dokchitser constants of an $\cR G$-lattice only depend on its genus.
\end{thm}
\begin{prop}
There exist at most 10 genera of indecomposable $\bZ D_{2p}$-lattices. Each genus has a
representative of the kind considered in Lemma \ref{lem:regConsts1}, Lemma \ref{lem:regConsts2} and
Lemma \ref{lem:regConsts3}.
\end{prop}
Our goal is to translate the Dokchitser constants into a certain index.
To that end, we now turn to the calculation of the index in $\Gamma$ of the submodule generated by
the various fixed submodules. The calculation is fairly similar to those of the Dokchitser
constants, but exhibits some new features. We will not give it in full detail but will give enough
examples to show the main techniques. The result is summarised in the following table,
where $C_2$ and $C_2'$ are two distinct subgroups of $D_{2p}$ isomorphic
to the cyclic group of order 2:
\begin{center}
\begin{tabular}{| r | l | l | l | l | l | l | l | l | l | l |}
\hline
$\Gamma$: & $\triv$ & $\epsilon$ & $\Delta$ & $A_i$ & $A_i'$ & $(A_i',\triv)$ & $(A_i,\epsilon)$ & $(A_i,\Delta)$ & $(A_i',\Delta)$ & $(A_i\oplus A_i',\Delta)$\\
\hline
$[\Gamma : \Gamma^{C_2}+ \Gamma^{C_2'}+ \Gamma^{C_p}]$: & 1 & 1 & 1 & $1 \;\forall i$ &
$p \;\forall i$ & $1 \;\forall i$ & $p \;\forall i$ & $p \;\forall i$ & $1 \;\forall i$ &
$p \;\forall i$\\
\hline
\end{tabular}
\end{center}

The assertions are clear for the first three modules in the list. For the others,
we begin by noting that the index is an invariant of the genus. More precisely the $l$-part of
the index for $\Gamma$ is equal to the $l$-part of the index for $\Gamma\otimes \bZ_{(l)}$ for any
prime $l$.
\begin{lem}
Let $\Gamma=A$. Then $\Gamma = \Gamma^{C_2}+ \Gamma^{C_2'}+ \Gamma^{C_p}$.
\end{lem}
\begin{proof}
We have already noted in the proof of Lemma \ref{lem:regConsts1} that $A$ is isomorphic
to the submodule of $\bZ[G/C_2] \otimes_\bZ \epsilon$ given by
$\left<1-\varpi^i: i = 1,\ldots,p-1\right>$ and that the submodule fixed by $\left<\sigma\right>$ is
\[\left<\varpi^i-\varpi^{p-i}: i = 1,\ldots,\frac{p-1}{2}\right>.
\]
It follows that the submodule
fixed by, say, $\left<\varpi^{-1}\sigma\varpi\right>$ is given by
\[\left<\varpi^{-1}(\varpi^i-\varpi^{p-i}): i = 1,\ldots,\frac{p-1}{2}\right>=
\left<\varpi^i-\varpi^{p-(i+2)}: i = 0,\ldots,\frac{p-3}{2}\right>.\]
These two are easily seen to
generate $A$. For example, by alternatingly summing elements from them one can obtain
\[1-\varpi = (1-\varpi^{-2}) + (\varpi^{-2}-\varpi^2) + (\varpi^{2}-\varpi^{-4}) +
\ldots +(\varpi^{p-3}-\varpi)\]
and similarly for all the other generators of $A$.
\end{proof}
The proof for $\Gamma=(A',\triv) = \bZ[G/C_2]$ is very similar in spirit and we will omit it.
\begin{lem}
Let $\Gamma=(A\otimes A',\Delta)\cong \bZ[G/1]$.
Then $[\Gamma : \Gamma^{C_2}+ \Gamma^{C_2'}+ \Gamma^{C_p}] = p$.
\end{lem}
\begin{proof}
The fixed submodules of $\bZ[G/1]$ under $\left<\sigma\right>$,
under $\left<\varpi\sigma\varpi^{-1}\right>$ and
under $\left<\varpi\right>$, respectively, are immediately seen to be
\[\left<\varpi^i+\sigma\varpi^i: i = 0,\ldots,p-1\right>,
\left<\varpi^i+\sigma\varpi^{i+2}: i = 0,\ldots,p-1\right>
\text{ and}\left<\sum_{i=0}^{p-1}\varpi^i,\sum_{i=0}^{p-1}\sigma\varpi^i\right>\]
and it is easy to check that
together these submodules generate the kernel of the surjective map
\[\bZ[G/1]\rightarrow \bZ/p \bZ,\;\;
\left(\sum_i \alpha_i \sigma\varpi^i + \sum_j \beta_j \varpi^j\right) \mapsto
\sum_i\alpha_i - \sum_j\beta_j \mod p.\]
This kernel is of index $p$ in $\bZ[G/1]$ and the claim is established.
\end{proof}
A similar proof, which we omit, shows the same for $\Gamma=A'$.

There are now several ways to finish the calculation.
For example, one can note that to compute the $p$-part of
the indices, we can localise everything at $p$ and use multiplicativity of indices in direct sums.
The $p$-parts of all the remaining indices then follow from the direct sum decompositions of
the proof of Lemma \ref{lem:regConsts3}.

We note that, by inspection, the quantity
\[\I(\Gamma) = \C_\Theta(\Gamma)\cdot[\Gamma : \Gamma^{C_2}+ \Gamma^{C_2'}+ \Gamma^{C_p}]^2\]
only depends on the rational representation $\Gamma\otimes \bQ$ and not on the lattice itself.
We deduce
\begin{lem}\label{lem:regConstIndex}
Write $\Gamma\otimes \bQ = \Lambda$.
Let $\triv$, $\sign$ and $\rho$ denote the irreducible rational
representations of $D_{2p}$, where $\rho$ is $(p-1)$-dimensional.\footnote{The abuse of notation in
using the same letters for the 1-dimensional rational representations and integral lattices in them
is very mild, since there is a unique integral lattice up to isomorphism in each of the two
representations.} Denote by $\left<\Lambda,.\right>$ the number of copies of a given irreducible
rational representation in $\Lambda$, analogous to the inner product of complex characters.
Then we have $\I(\triv) = 1/p$, $\I(\sign)=p$, $\I(A_i) = \I(A_i') = p$ for all $i$, and for any $\Gamma$,
we have $\I(\Gamma) = p^{\left<\Lambda, \sign\right>+\left<\Lambda, \rho\right>-\left<\Lambda, \triv\right>}$.
\end{lem}
This identity will be crucial in proving Theorem \ref{thm:unitIndex}.
\section{Class number relations - main results}\label{sec:mainRes}
In this section we will collect the results obtained so far to prove the main theorems.
\subsection{Possible values of regulator quotients}
We will begin by establishing Theorem \ref{thm:noP}.
\begin{thm}\label{thm:noPInLambda}
Let $G$ be a finite group, let $N$ be a normal subgroup such that $G/N$ is cyclic, let $l$ be a prime
number not dividing the order of $N$. Let $F/K$ be a Galois extension of number fields with Galois
group $G$ and $\Theta = \sum_i H_i - \sum_j H_j'$ a $G$-relation. Let $S$ be a finite $G$-stable
set of places of $F$ including all the Archimedean ones. Recall the notation
\[
\lambda(H)=\#\ker\left(\sH^1(H,\mu(F))\rightarrow \sH^1(H,\Sunits{F})\right)
\]
for $H\leq G$. Then
\[\ord_l\left(\prod_i \lambda(H_i)\big/\prod_j \lambda(H_j')\right) = 0.\]
\end{thm}
\begin{proof}
For a subgroup $H$ of $G$, define $\tau_F(H)$ as
\[\tau_F(H) = \ker\left(\sH^1(H,\mu(F))\rightarrow \sH^1(H,\Sunits{F})\right),\]
so that $\lambda(H) = \#\tau_F(H)$.
The inflation-restriction exact sequence gives us the following commutative diagram with exact rows:
\[
\xymatrix{
0 \ar[r] & \sH^1(H/(H\cap N),\mu(F^{H\cap N})) \ar[d]\ar[r] & \sH^1(H,\mu(F)) \ar[d]\ar[r] & \sH^1(H\cap N,\mu(F))\ar[d]\\
0 \ar[r] & \sH^1(H/(H\cap N),\Sunits{F^{H\cap N}})       \ar[r] & \sH^1(H,\Sunits{F})       \ar[r] & \sH^1(H\cap N,\Sunits{F})
}
\]
where the commutativity is obvious on the level of co-cycles. Hence we get the exact sequence
\[0\rightarrow \tau_{F^N}(HN/N) \rightarrow \tau_F(H) \rightarrow \tau_F(H\cap N).\]
But the $l$-part of $\lambda(H\cap N)$ is trivial, since $l$ does not divide $|H\cap N|$ and so we
see that $\lambda(H)[l^\infty]= \lambda(HN/N)[l^\infty]=\#\tau_{F^N}(HN/N)[l^\infty]$.
Since $G/N$ is cyclic and therefore has
no non-trivial relations, and by applying \cite[Theorem 2.36 (q)]{TVD-2} with
$\phi_{G/N}(HN/N)=\lambda(HN/N)[l^\infty]$, $\lambda(H)[l^\infty]$ vanishes in
relations and we are done.
\end{proof}
\begin{cor}
Under the hypotheses of the theorem, we have
\[\ord_l\left(\prod_i R_S(F^{H_i})\big/\prod_j R_S(F^{H_j'})\right) = 0.\]
\end{cor}
\begin{proof}
This is an immediate consequence of Proposition \ref{prop:newReg}, Proposition
\ref{prop:normalSgpCycQuo} and the above theorem.
\end{proof}
\begin{cor}
Under the hypotheses of the theorem, we have an equality of the $l$-parts of class numbers:
\[\prod_i h_S(F^{H_i})_l = \prod_j h_S(F^{H_j'})_l.\]
\end{cor}
\begin{proof}
We only need to establish that the $l$-part of the quotient
$\prod_i w(F^{H_i})/\prod_j w(F^{H_j'})$ is trivial. If $l\neq 2$, then this is true in general as
observed by Brauer \cite[\S 2]{Bra-51}. If $l=2$, then $w(F^H)=w(F^{HN})$ and the latter vanishes
in relations by exactly the same argument as in the proof of Theorem \ref{thm:noPInLambda}.
\end{proof}
\begin{rmrk}
One could also deduce both Corollaries directly, without using Theorem \ref{thm:noPInLambda}, from
the work of Boltje and of Bley and Boltje on Mackey functors (see e.g. \cite[Corollary 2.4]{Bol-97}
or \cite{BB-04}) combined with the proof of Proposition \ref{prop:normalSgpCycQuo}.
\end{rmrk}
\subsection{Unit index formula for $D_{2p}$-extensions}
We will first prove Theorem \ref{thm:unitIndex} and then
explain how to deduce a formula for $D_{2q}$ where $q$ is any odd integer. We will not actually
write down the formula for
$D_{2q}$ because it is less enlightening when it is written out than its conceptual idea. The
interested reader should have no difficulties in writing it down for any specific case.
Let $F/k$ be a Galois extension of number fields with Galois group
$G=D_{2p}$ for $p$ an odd prime.
Let $K$ be the intermediate quadratic extension and $L\neq L'$ two
intermediate extensions of degree $p$ over $k$. As in the previous section,
denote by $ \varpi$ an element of order $p$ in $G$
and let $\sigma$ be the involution that fixes $L$. Let $\Theta$ be the relation from Example
\ref{ex:d2p} and let $S$ be a $G$-stable set of primes of $F$. Our main tool
is the compatibility statement between Artin formalism and the analytic class number formula given
by equation (\ref{eq:classNoFormula}). We will first show that in our case,
$\frac{w(F)w(k)^2}{w(K)w(L)^2}=1$. Indeed, since the extension $L/k$ is not Galois, it can not be
obtained by adjoining roots of unity. Since it has no intermediate extensions, we see that $w(L)=w(k)$.
Similarly, if $F$ was obtained from $K$ by adjoining roots of unity, then adjoining these same roots
to $k$ would give an extension of degree $p$ or $2p$. But the former is not possible by what we have
just said and the latter would imply that $F/k$ is abelian. So $w(F)=w(K)$ and our claim follows.

Set $\Gamma$ to be the Galois module $U_S(F)$ given by the $S$-units of $F$ modulo torsion and write
$\Lambda=\Gamma\otimes \bQ$. We will now invoke Proposition \ref{prop:newReg}. Note that in our
situation, the only subgroup of $G$ for which
$\C_\Theta(\bZ[G/H])\neq 1$ is $G$ itself and that $\C_\Theta(\bZ[G/G]) = \C_\Theta(\triv) = 1/p$.
So $\prod_{\fp\in S|_k}\C_\Theta(\bZ[G/D_\fp]) = p^{-\#\{\fp\in S|_k:\;D_\fp = G\}}$. Set
$a(F/k,S)=\#\{\fp\in S|_k:\;D_\fp = G\}$.
Then we have, using the notation from Proposition \ref{prop:newReg} and from Lemma
\ref{lem:regConstIndex},
\begin{eqnarray}
\frac{h_S(F)h_S(k)^2}{h_S(K)h_S(L)^2} & \stackrel{\text{eqn. }(\ref{eq:classNoFormula})}{=} &
\frac{R_S(K)R_S(L)^2}{R_S(F)R_S(k)^2}\nonumber\\
& \stackrel{\text{Prop. }\ref{prop:newReg}}{=} &
\left(\C_\Theta(\triv)\big/\left(\C_\Theta(\Gamma)\prod_{\fp\in S|_k}\C_\Theta(\bZ[G/D_\fp])
\right)\right)^{1/2}\cdot\frac{\lambda(C_p)\lambda(C_2)^2}{\lambda(1)\lambda(G)^2}\nonumber\\
& = & \left(p^{a(F/k,S)-1}\big/\C_\Theta(\Gamma)\right)^{1/2}\cdot
\frac{\lambda(C_p)\lambda(C_2)^2}{\lambda(1)\lambda(G)^2}\nonumber\\
& \stackrel{\text{Lem. }\ref{lem:regConstIndex}}{=} & \left(p^{a(F/k,S)-1}\cdot [\Gamma: \Gamma^{C_2}+\Gamma^{C_2'}+\Gamma^{C_p}]^2
\big/p^{\left<\Lambda, \sign\right>+\left<\Lambda, \rho\right>-\left<\Lambda, \triv\right>}\right)^{1/2}\times\nonumber\\
& & \frac{\lambda(C_p)\lambda(C_2)^2}{\lambda(1)\lambda(G)^2}\nonumber\\
& = & \left(p^{r_S(k)-(r_S(K)-r_S(k))-(r_S(F)-r_S(K))/(p-1)+a(F/k,S)-1}\right)^{1/2}\times\nonumber\\
& & [\Gamma: \Gamma^{C_2}+\Gamma^{C_2'}+\Gamma^{C_p}]\cdot\frac{\lambda(C_p)\lambda(C_2)^2}{\lambda(1)\lambda(G)^2}.\label{eq:mainFormula}
\end{eqnarray}

Recall that $\lambda(H) = [U_S(F)^H:U_S(F^H)] = \#\ker(\sH^1(H,\mu(F))\rightarrow \sH^1(H,\Sunits{F}))$.
As we have discussed before Lemma \ref{lem:index}, this is trivial for all $H\leq G$ if neither
$F/L$ nor $F/K$ is obtained by adjoining a root of a non-torsion $S$-unit.

It remains to compute $\lambda(H)$ for all
$H\leq G$ and to compare the $\Gamma$-index in equation (\ref{eq:mainFormula}) with the actual
unit index appearing in Theorem \ref{thm:unitIndex}. Since the roots
of unity $\mu(F)$ are contained in $\Sunits{K}$ as remarked above, we have
\begin{eqnarray}
  \lefteqn{[\Sunits{F} : \Sunits{L}\Sunits{L'}\Sunits{K}] = }\nonumber\\
  & = &[\Sunits{F}/\mu(F) : \Sunits{L}\Sunits{L'}\Sunits{K}/\mu(F)]\nonumber\\
  & = & \left[\Sunits{F}/\mu(F):
\left(\mu(F)\Sunits{L}/\mu(F)\right)
\left(\mu(F)\Sunits{L'}/\mu(F)\right)\left(\mu(F)\Sunits{K}/\mu(F)\right)\right]\nonumber\\
  & = & \left[\Sunits{F}/\mu(F):
\left(\Sunits{L}/\mu(F)\cap \Sunits{L}\right)
\left(\Sunits{L'}/\mu(F)\cap \Sunits{L'}\right)
\left(\Sunits{K}/\mu(F)\cap \Sunits{K}\right)\right]\nonumber\\
  & = & \left[\Sunits{F}/\mu(F):
\left(\Sunits{L}/\mu(L)\right)\left(\Sunits{L'}/\mu(L')\right)\left(\Sunits{K}/\mu(K)\right)\right]\\
  & = & \left[U_S(F):U_S(F^{C_2})U_S(F^{C_2'})U_S(F^{C_p})\right],
\label{eq:unitIndexComp}
\end{eqnarray}
where for $H\leq G$, $\Sunits{F^H}/\mu(F^H)$ is identified with its image in $\Sunits{F}/\mu(F)$
under the obvious inclusion map. Thus, to compare the $\Gamma$-index with the unit index, we need
to compute
\[i(F/k,S) =[U_S(F)^{C_2}U_S(F)^{C_2'}U_S(F)^{C_p}:U_S(F^{C_2})U_S(F^{C_2'})U_S(F^{C_p})].\]
We will consider various different cases. The remaining computations are
summarised in the following two lemmata:
\begin{lem}
We have $\lambda(1) = 1$; $\lambda(C_2)\in \{1,2\}$; $\lambda(C_p)\in \{1,p\}$
with $\lambda(C_p)=p$ if and only if $F=K(\sqrt[p]{u})$ for a non-torsion $S$-unit
$u\in \Sunits{K}$;
$\lambda(G)\in\{1,2,p,2p\}$ with $\lambda(G)$ divisible by 2 if and only if $\lambda(C_2)=2$ and
divisible by $p$ if and only if $L=k(\sqrt[p]{u'})$ for a non-torsion $S$-unit $u'\in \Sunits{k}$.
\end{lem}
\begin{proof}
It is clear that $\lambda(1) = 1$. Since, for any $H\leq G$, $\lambda(H)$ is the order
of a subgroup of $\sH^1(H,\mu(F))$ and since this cohomology group is cyclic and annihilated
by $|H|$, we deduce that $\lambda(H)$ divides $|H|$ for $H\in \{C_2,C_p,G\}$.
The 2-part of the $\lambda$-quotient vanishes by Theorem \ref{thm:noPInLambda}.

By definition, an element of $U_S(F)^{C_p}/U_S(F^{C_p})$ of order $p$ is represented by a
non-torsion $S$-unit
$v\in \Sunits{F}\backslash \mu(F)\Sunits{K}$ such that $\varpi(v) = \zeta v$ and
$v^p = \bar{\zeta}x$ for $\zeta,\bar{\zeta}\in \mu(F)$ and $x\in K$. But
$\mu(F)\subset \Sunits{K}$, so these conditions are equivalent to
$v\in \Sunits{F}\backslash \Sunits{K},\;v^p\in K$. Thus, $\lambda(C_p)=p$ if and only if
$F=K(\sqrt[p]{u})$, where $u$ is a non-torsion $S$-unit in $\Sunits{K}$.

Also, an element of $U_S(F)^{G}/U_S(F^{G})$ of order $p$ is represented by an $S$-unit
$v'\in \Sunits{F}\backslash \mu(F)\Sunits{k}$ such that $\varpi(v') = \zeta_1 v'$,
$\sigma(v') = \zeta_2 v'$ and $v'^p = \bar{\zeta}x$
for $\zeta_1,\zeta_2,\bar{\zeta}\in \mu(F)$ and $x\in \Sunits{k}$. If $L=k(\sqrt[p]{u'})$
for a non-torsion $S$-unit $u'\in \Sunits{k}$, then the conditions are satisfied for
$v'=\sqrt[p]{u'}$, so in this case $\lambda(G)$ is divisible by $p$. Conversely, let
$v'\in \Sunits{F}\backslash \mu(F)\Sunits{k}$ represent an element of order $p$ and let
$\zeta_1,\zeta_2,\bar{\zeta}$ and $x$ be as above. We need to find
$\tilde{v}\in \Sunits{F}\backslash \mu(F)\Sunits{k}$ satisfying the same conditions, but with
$\tilde{v}^p\in \Sunits{k}$ (and not merely in $\mu(F)\Sunits{k}$). Consider
$\tilde{v} = \Norm_{F/L}(v') = v'\sigma(v') = \zeta_2v'^2.$ Clearly, it is fixed by $G$ up to roots
of unity, since $v'$ is, and also $v'^p\in \mu(F)\Sunits{k} \subseteq \Sunits{K}$ implies that
$\tilde{v}^p\in \Sunits{k}$. So we only need to show that $\tilde{v}\notin \mu(F)\Sunits{k}$.
But if $\tilde{v}\in \mu(F)\Sunits{k}\subseteq\Sunits{K}$, then in fact $\zeta_2v'^2=
\tilde{v}\in K\cap L = k$, so $v'^2\in \mu(F)\Sunits{k}$, contradicting the assumption that $v'$
represents an element of order $p$ in $U_S(F)^{G}/U_S(F^{G})$.
\end{proof}

\begin{lem}
The index $i(F/k,S)=[U_S(F)^{C_2}U_S(F)^{C_2'}U_S(F)^{C_p}:U_S(F^{C_2})U_S(F^{C_2'})U_S(F^{C_p})]$
is equal
to $p$ if $F=K(\sqrt[p]{u})$ for a non-torsion $S$-unit $u\in \Sunits{K}$ and $L$ is not
obtained by adjoining a non-torsion $S$-unit to $k$, and is 1 otherwise.
\end{lem}
\begin{proof}
The statement clearly holds if $\lambda(H)=1$ for all $H\leq G$.
Next, by the previous lemma, any non-trivial element in $U_S(F)^{C_2}/U_S(F^{C_2})$ can be
represented by an element of $U_S(F)^G$. In particular, this representative is fixed by $C_p$ up
to roots of unity, so gives an element of $U_S(F)^{C_p}/U_S(F^{C_p})$. Since this latter group has
no 2-torsion, we deduce that any non-trivial class in $U_S(F)^{C_2}/U_S(F^{C_2})$ is represented
by an element of $U_S(F^{C_p})$ and so the index $i(F/k,S)$ is never divisible by 2.

Similarly, by the same lemma, if $L=k(\sqrt[p]{u'})$ for a non-torsion $S$-unit $u'\in \Sunits{k}$,
then a generator
of $U_S(F)^{C_p}/U_S(F^{C_p})$ can be represented by an element of $U_S(F)^{G}$, which then gives
an element of $U_S(F)^{C_2}/U_S(F^{C_2})$. This group has no $p$-torsion, so the generator
of $U_S(F)^{C_p}/U_S(F^{C_p})$ is represented by an element of $U_S(F^{C_2})$ and the index
$i(F/k,S)$ is not divisible by $p$ in this case, hence trivial.

Finally, suppose that
$F=K(\sqrt[p]{u})$ for a non-torsion $S$-unit $u\in \Sunits{K}$, but that $L$ is not obtained from
$k$ in this way. We will show that then, $\sqrt[p]{u}$ represents a non-trivial coset of
\[U_S(F)^{C_2}U_S(F)^{C_2'}U_S(F)^{C_p}/U_S(F^{C_2})U_S(F^{C_2'})U_S(F^{C_p}),\]
necessarily of order $p$ in this quotient.
Assume for a contradiction that $\sqrt[p]{u}=u_Lu_{L'}u_K$, where
$u_M \in \Sunits{M}$ for $M=L,L',K$. Recall, that $\varpi$ denotes an element of $G$ of order $p$.
Let $\zeta_p$ be a primitive $p$-th root of unity in $K$, which must exist since otherwise
$F/K$ would not be Galois. Since
$\zeta_pu_Lu_{L'}u_K=\zeta_p\sqrt[p]{u}=\varpi(\sqrt[p]{u}) =u_K\varpi(u_Lu_{L'})$, we may replace
$\sqrt[p]{u}$ by $\sqrt[p]{u}/u_K$ and assume without loss of generality that $u_K=1$. Consider
the images $w$, $w_L$ and $w_{L'}$ of $\sqrt[p]{u}$, $u_L$ and $u_{L'}$, respectively,
in $\Sunits{F}\otimes_\bZ \bQ$, where $\Sunits{F}$ is regarded
as a $\bZ$-module by virtue of being a finitely generated abelian group. Recall the notation
$\triv$, $\sign$ and $\rho$ of Lemma \ref{lem:regConstIndex} for the irreducible rational
representations of $D_{2p}$.
Since $w$ is fixed by $C_p$ and since the $C_p$-invariant subspace of $\rho$ is trivial (this is
true for the complex irreducible two-dimensional representations, of which $\rho$ is the sum),
the projection of $w$ onto the $\rho$-isotypical component of $\Sunits{F}\otimes_\bZ \bQ$ is trivial.
Since the $C_2$-invariant subspace and the $C_2'$-invariant subspace of $\rho$ are linearly
independent (as can again be seen on the level of the complex irreducible summands of $\rho$),
the projections of $w_L$ and of $w_{L'}$ onto the $\rho$-isotypical component must also be trivial.
But also, the $C_2$-invariant subspace of the $\sign$-isotypical component is trivial, so we deduce
that $w_{L}$ is in the $\triv$-component, in other words that $G$ acts on $u_L$ by multiplying
by roots of unity. So either $u_L\in k$ and thus $u\in L'$,
or $L$ is obtained from $k$ by obtaining the $p$-th
root of a non-torsion $S$-unit, both possibilities contradicting the assumptions.
\end{proof}
Combining the two lemmata with equations (\ref{eq:mainFormula}) and (\ref{eq:unitIndexComp})
completes the proof of Theorem \ref{thm:unitIndex}.

\subsection{A formula for $D_{2q}$ for $q$ any odd integer}
Throughout this subsection we fix the following notation:
\begin{notation*}
In this subsection we will drop the subscript $S$ from $\Sunits{M}$ and write $\Units{M}$ instead.
The set-up is as follows

\begin{tabular}{ l l }
$q$ & $\prod_{i=1}^n p_i$ for $p_i$ odd primes, not necessarily distinct;\\
$G$ & dihedral group with $2q$ elements, $D_{2q}=\left<a,b\;|\; a^q=b^2=(ab)^2=1\right>;$\\
$F/k$ & a Galois extension of number fields with Galois group $G$;\\
$K$ & $F^{\left<a\right>}$;\\
$L$ & $F^{\left<b\right>}$;\\
$L'$ & $F^{\left<ab\right>}$;\\
$S$ & Galois stable set of places of $F$ including the Archimedean ones;\\
\end{tabular}\vspace{\baselineskip}

For each $j\in \{0,\ldots,n\}$ define

\begin{tabular}{ l l }
$\varepsilon_j$ & $\prod_{i=1}^j p_i$ (1 if $j=0$);\\
$C_j$  & $\left<a^{\varepsilon_j}\right>$, the unique subgroup of $\left<a\right>$ of index $\prod_{i=1}^jp_i$ (1 if $j=0$);\\
$D_j$, $D_j'$ & the dihedral groups generated by $C_j$ and $b$ and by $C_j'$ and $b$,\\
& respectively.\\
\end{tabular}
\end{notation*}
With this notation, $F^{C_j}/F^{D_{j-1}}$ is an intermediate Galois extension with Galois group
$D_{2p_j}$ for $j\in \{1,\ldots,n\}$ and so Theorem \ref{thm:unitIndex} applies to this extension.
By taking the product of the unit index formula over $j=1,\ldots,n$, we obtain that
\begin{eqnarray}\label{eq:product}
\frac{h_S(F)h_S(k)^2}{h_S(K)h_S(L)^2} = \prod_{j=1}^n \left(p_j^{\alpha_j}\times
\left[\Units{F^{C_j}}:\Units{F^{C_{j-1}}}\Units{F^{D_j}}\Units{F^{D_j'}}\right]\right),
\end{eqnarray}
where $\alpha_j$ are the corresponding exponents of $p_j$ from Theorem \ref{thm:unitIndex}.

Before investigating the unit index, we will give a more conceptual explanation of this formula.
We have the $G$-relation
\[\Theta=1 - 2C_2 - C_{p^n} + 2G.\]
As in the case of $D_{2p}$, the corresponding quotient of numbers of roots of unity
$\frac{w(F)w(k)^2}{w(K)w(L)^2}$ is trivial, because if $F$ contains a root of unity, then adjoining
this root to $k$ gives an abelian Galois extension of $k$ which must therefore be contained in $K$.
Thus $w(F)=w(K)$ and $w(L)=w(k)$. So, using equation
(\ref{eq:classNoFormula}), we see that 
\[\frac{h_S(F)h_S(k)^2}{h_S(K)h_S(L)^2} = \frac{R_S(K)R_S(L)^2}{R_S(F)R_S(k)^2}\]
and Proposition \ref{prop:newReg} implies that
\[\frac{h_S(F)h_S(k)^2}{h_S(K)h_S(L)^2} =
\left(\frac{\C_\Theta(\triv)}{\C_\Theta(\Gamma)\cdot\prod_{\fp\in S|_k}\C_\Theta(\bZ[G/D_\fp])}
\right)^{1/2}\frac{\lambda(C_{p^n})\lambda(C_2)^2}{\lambda(1)\lambda(G)^2}.\]
This time, we do not have a classification of all
indecomposable integral representations of $G$ at our disposal (in fact the number of their isomorphism
classes is infinite when $q$ is not cube-free). Instead, to replace the Dokchitser constant by a unit
index, we break up the Dokchitser constant into Dokchitser constants of $D_{2p_i}$-representations and then use Lemma \ref{lem:regConstIndex}.
We begin by an obvious Lemma:
\begin{lem}\label{lem:onlyFixedSubmod}
Let $G$ be any finite group, $\Theta = \sum_{i\in I} n_iH_i$ any $G$-relation with $n_i$ non-zero integers
and $\Gamma$ any $\cR$-free $\cR G$-module. Set $H=\cap_{i \in I}H_i$. Then
\[\C_\Theta(\Gamma) = \C_\Theta(\Gamma^H).\]
\end{lem}
\begin{proof}
This is clear from the definition of Dokchitser constants, since elements of $\Gamma$ that are not fixed
by any of the subgroups occurring in the relation do not contribute to the Dokchitser constant.
\end{proof}
For each integer $j\in \{1,\ldots,n\}$ we have the $G$-relation
$\Theta_j = C_j -2D_j - C_{j-1} + 2D_{j-1}$.
We see immediately that $\Theta = \sum_{j=1}^n \Theta_j$ and so by Proposition \ref{prop:multiplicativity}
we have
\[\C_\Theta(\Gamma) = \prod_{j=1}^n \C_{\Theta_j}(\Gamma).\]
For each $j$, $\Theta_j$ is induced from the corresponding relation in $D_{j-1}$ and so by Proposition
\ref{prop:restrictionInduction} we have $\C_{\Theta_j}(\Gamma) = \C_{\Theta_j}(\Gamma\!\!\downarrow_{D_{j-1}})$,
where on the right hand side $\Theta_j$ is viewed as a $D_{j-1}$-relation. Moreover, by Lemma \ref{lem:onlyFixedSubmod}
we have $\C_{\Theta_j}(\Gamma) = \C_{\Theta_j}((\Gamma\!\!\downarrow_{D_{j-1}})^{C_j})$. Now,
$(\Gamma\!\!\downarrow_{D_{j-1}})^{C_j}$ can be considered as a
$(D_{j-1}/C_j\cong D_{2p_j})$-module $\bar{\Gamma}_j$, and since $\Theta_j$ is in fact lifted
from the $D_{2p_j}$-quotient of $D_{j-1}$, we have from Proposition \ref{prop:restrictionInduction}
\[\C_\Theta(\Gamma) = \prod_{j=1}^n \C_{\Theta_j}(\bar{\Gamma}_j),\]
where each factor is now a Dokchitser constant in $D_{2p_j}$. Applying Lemma \ref{lem:regConstIndex}
and the discussion of the $D_{2p}$ case recovers equation (\ref{eq:product}).

Ideally, we would like to replace the product of the unit indices by the index
\[[\Units{F}:\Units{K}\Units{L}\Units{L'}].\]
However, the right hand side of equation (\ref{eq:product})
depends on more than this one index and some correction terms will be necessary.

Write $K_j=F^{C_j}$, $L_j=F^{D_j}$ and $L_j'=F^{D_j'}$, so that for example $K_0=K$, $K_n=F$, $L_0=k$
and $L_n=L$:
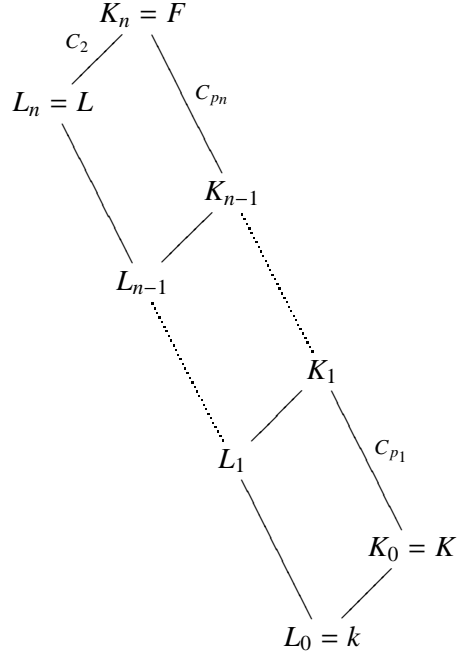
\[
\xymatrix@!C@!=0.8pc{
& K_n=F \ar@{-}[ddr]^{C_{p_n}} \ar@{-}[dl]_{C_2}\\
L_n=L \ar@{-}[ddr]\\
& & K_{n-1} \ar@{.}[ddr] \ar@{-}[dl]\\
& L_{n-1} \ar@{.}[ddr]\\
& & & K_1 \ar@{-}[ddr]^{C_{p_1}} \ar@{-}[dl]\\
& & L_1 \ar@{-}[ddr]\\
& & & & K_0=K \ar@{-}[dl]\\
& & & L_0=k\\
}
\]
First, note that for any group $X$ and
any normal subgroups $Y$ and $Z$, we have $|X/Y|=|YZ/Y|\cdot|X/YZ|=|Z/(Y\cap Z)|\cdot|X/Z\big/YZ/Z|$,
provided all the quotients are finite. Applying this in step $\dagger$ below with
$X=\Units{L_n}\Units{L_n'}\Units{K_1}$,
$Y=\Units{L_n}\Units{L_n'}\Units{K_0}$ and $Z=\Units{K_1}$, we get
\begin{eqnarray*}
[\Units{F}:\Units{L_n}\Units{L_n'}\Units{K_0}] & = &[\Units{F}:\Units{L_n}\Units{L_n'}\Units{K_1}]
\times \;[\Units{L_n}\Units{L_n'}\Units{K_1}:\Units{L_n}\Units{L_n'}\Units{K_0}]\\
& \stackrel{\dagger}{=} & [\Units{F}:\Units{L_n}\Units{L_n'}\Units{K_1}]
\times \;[\Units{L_n}\Units{L_n'}\Units{K_1}\cap \Units{K_1}:\Units{L_n}\Units{L_n'}\Units{K_0}\cap \Units{K_1}]\times\\
& & \times \;[\Units{L_n}\Units{L_n'}\Units{K_1}/\Units{K_1}:\Units{L_n}\Units{L_n'}\Units{K_0}\Units{K_1}/\Units{K_1}]\\
& = & [\Units{F}:\Units{L_n}\Units{L_n'}\Units{K_1}]
\times \;[\Units{K_1}:\Units{L_n}\Units{L_n'}\Units{K_0}\cap \Units{K_1}]\\
& = & [\Units{F}:\Units{L_n}\Units{L_n'}\Units{K_1}]\times\\
& & \times \;[\Units{K_1}:\Units{L_1}\Units{L_1'}\Units{K_0}]
\big/\; [\Units{L_n}\Units{L_n'}\Units{K_0}\cap \Units{K_1}: \Units{L_1}\Units{L_1'}\Units{K_0}].
\end{eqnarray*}

Repeating this inductively yields

\begin{eqnarray*}
\lefteqn{\prod_{j=1}^n[\Units{K_j}:\Units{L_j}\Units{L_j'}\Units{K_{j-1}}]=}\\
& & [\Units{F}:\Units{L_n}\Units{L_n'}\Units{K_0}]
\times\;\prod_j[\Units{L_n}\Units{L_n'}\Units{K_{j-1}}\cap \Units{K_{j}}: \Units{L_j}\Units{L_j'}\Units{K_{j-1}}].
\end{eqnarray*}
Finally, substituting this in equation (\ref{eq:product}) gives the sought for unit index formula.

\section{Examples}\label{sec:examples}
We first derive some easy consequences of Theorem \ref{thm:unitIndex}:
\begin{cor}
Let $F/\bQ$ be a Galois extension with Galois group $D_{2p}$ for $p$ an odd prime. Let $K$ be the
quadratic subfield and let $L$ and $L'$ be distinct intermediate extensions of degree $p$ over $\bQ$. Let $r(K)$ be the rank
of the units in $K$, which is either 0 or 1. Then we have
\[\frac{h(F)p^{r(K)+1}}{h(K)h(L)^2}=[\Units{F}:\Units{K}\Units{L}\Units{L'}].\]
\end{cor}
This is the formula derived by Halter-Koch in \cite{HK-77}.
\begin{proof}
In this case, $S$ consists of the Archimedean primes and none of them have decomposition group $D_{2p}$.
Moreover, $F/K$ cannot be obtained by adjoining a $p$-th root of a fundamental unit, since for that
$K$ has to contain the $p$-th roots of unity and have unit rank 1, which is impossible. Finally,
$r(F)=p(r(K)+1)-1 = p\cdot r(K)+p-1$. So formula (\ref{eq:mainFormula}) simplifies to the stated form.
\end{proof}
\begin{cor}
Let $F/k$ be a Galois extension of number fields with Galois group $G=D_{2p}$ for $p$ an odd prime,
let $K$ be the intermediate quadratic extension and let $L$ and $L'$ be distinct intermediate extensions of degree $p$. Let
$S$ be a $G$-stable set of primes of $F$ including the Archimedean ones such that their decomposition
groups do not contain $C_p$. Also assume that $F/K$ is not obtained by adjoining a $p$-th root of
a non-torsion $S$-unit. Then
\[
\frac{h_S(F)h_S(k)^2p^{r_S(K)+1-r_S(k)}}{h_S(K)h_S(L)^2} = [\Sunits{F}:\Sunits{L}\Sunits{L'}\Sunits{K}].
\]
\end{cor}
The condition that $F\neq K(\sqrt[p]{u})$ for a non-torsion $S$-unit $u$ of $K$ is for example
satisfied when $K$ does not contain the $p$-th roots of unity
or when $F/K$ is unramified at $p$. In particular, the corollary applies when $F/K$ is unramified, so
this includes the case considered by Lemmermeyer in \cite[Theorem 2.2]{Lem-05}.
\begin{proof}
We again have that $r_S(F) = p\cdot r_S(K) + p-1$, since all the places in $S$ are assumed to split in
$F/K$, and the claim is a direct consequence of formula (\ref{eq:mainFormula}).
\end{proof}
In particular cases we can use the classification of integral representations of $D_{2p}$ to express
the Galois structure of the units modulo torsion in terms of the class number quotient. This
has been explored when the base field is $\bQ$ and $S$ contains only the Archimedean place,
e.g. in \cite{Mos-75}. We will give some more examples in the more general setting.
\begin{ex}
Let $k$ be a real quadratic field and let $F/k$ be a Galois extension with Galois group $G=D_{2p}$.
As before, let $K$ be the intermediate quadratic extension and let $L$ be an intermediate extension
of degree $p$ and take $\Gamma$ to be the integral $G$-representation given by the units of $F$ modulo
torsion (or more precisely their usual logarithmic embedding into $\bR^{r(F)+1}$). Assume that $F/K$ is
not obtained by adjoining a $p$-th root of a non-torsion unit of $K$. Further, assume for simplicity
that $K$ is totally complex. Then, $r(k)=r(K)=1$ and $r(F)=2p-1$. So the $\bQ G$-representation
given by $\Gamma\otimes \bQ$ contains one copy of the trivial representation and two copies of the
$p-1$ dimensional irreducible representation. Using the notation from section \ref{sec:dihGroups} we
have the following possible $\bZ G$-module structures for $\Gamma$ together with the corresponding
class number quotients:
\[\begin{array}{l l}
\Gamma & \frac{h(F)h(k)^2}{h(K)h(L)^2}\\
\hline\\
A_i\oplus A_i\oplus \triv & 1/p\\
A_i\oplus A_i'\oplus \triv & 1\\
A_i'\oplus A_i'\oplus \triv & p\\
A_i\oplus (A_i',\triv) & 1/p\\
A_i'\oplus (A_i',\triv) & 1
\end{array}\]
where the values of the class number quotients follow from Proposition \ref{prop:newReg} and the
computation of Dokchitser constants in section \ref{sec:dihGroups}. In particular, we see that if
the class number quotient is $p$, then this determines the genus of the integral representation $\Gamma$.
We remind the reader that by the classification of integral representations, the number of the representations
$A_i$ in the same genus is equal to the class number of $\bQ(\zeta_p)^+$. This is known to be 1 for
$p\leq 67$ and conjectured to be 1 for $p\leq 157$ (this conjecture is implied by the generalised
Riemann hypothesis), so for 'small' $p$ the class number quotient can
sometimes completely determine the Galois module structure of the units modulo torsion.

If $K$ is not totally complex, then the same kind of analysis applies but the rank of the units of $F$
is larger and there are more possibilities to consider.
\end{ex}
\begin{ex}
In the previous example we have seen how, using our general result, we can apply Moser's reasoning from
\cite{Mos-75} to base fields, different from $\bQ$. We will now show how the generalisation to $S$-units can
be useful to complement Moser's results. Let $F/\bQ$ be a $D_{2p}$-extension with $K$, $L$ and
$\Gamma$ as above. If $K$ is imaginary, then $r(K)=0$ and $\Gamma\otimes \bQ$ only contains one copy
of the irreducible $(p-1)$-dimensional representation. By the classification of integral representations
and the computation of their Dokchitser constants in section \ref{sec:dihGroups}, we see that the
class number quotient is either $1$ or $1/p$ and in either case it determines the genus of $\Gamma$.
However, when $K$ is real, we have the following possibilities for $\Gamma$ together with the
corresponding class number quotients:
\[\begin{array}{l l l}
\text{number} & \Gamma & \frac{h(F)h(k)^2}{h(K)h(L)^2}\\
\hline\\
(1) & A_i\oplus A_i\oplus \sign & 1/p^2\\
(2) & A_i\oplus A_i'\oplus \sign & 1/p\\
(3) & A_i'\oplus A_i'\oplus \sign & 1\\
(4) & A_i\oplus (A_i,\sign) & 1/p\\
(5) & A_i'\oplus (A_i,\sign) & 1
\end{array}\]
We see that if the class number quotient is $1/p^2$, then the genus of $\Gamma$ is again determined
(and therefore the whole Galois module structure of $\Gamma$ is determined if $p\leq 67$, as remarked
in the previous example). However, if the class number quotient is 1 or $1/p$, then we are left with two
possibilities. But sometimes, looking at $S$-class numbers can resolve the ambiguity. Let $q$ be a
prime number which is inert or ramified in $K/\bQ$ and ramified in $F/K$. Let $S$ consist of the infinite
places of $F$ and the places above $q$. Let $\Gamma_S$ be the Galois module given by the $S$-units of
$F$ modulo torsion. Then $\Gamma_S\otimes \bQ$ contains one copy of the trivial representation, one
copy of the non-trivial one-dimensional representation and two copies of the irreducible
$(p-1)$-dimensional representation. Also, $\Gamma_S$ contains $\Gamma$ as a saturated sublattice and
the possible Galois module structures of $\Gamma_S$ restrict the possibilities for
$\Gamma$. For example, if the $S$-class number quotient is $1/p$, then writing out the list of possibilities
for $\Gamma_S$ (there are 16) we see that $\Gamma$ is given either by number (1) or by (2) and the
two have different
 class number quotients. Here is a concrete example: let $F$ be the splitting field
of the irreducible cubic polynomial
\[f(x) = x^3-34x-6.\]
The Galois group of $F/\bQ$ is $S_3$ and the
class number quotient is $1/3$. Thus, the Galois module structure of the units of $F$ modulo roots of
unity is either (2) or (4) from the above list. Now, let $S$ consist of the infinite places of $F$
and the unique place above 2. Then, the $S$-class number quotient is also $1/3$ and so the Galois
module structure of the units of $F$ modulo the roots of unity must be the one numbered (2) in the 
list.
\end{ex}
\addcontentsline{toc}{section}{Bibliography}
\bibliographystyle{classnorels}
\bibliography{ClassNoRelations}

\begin{thebibliography}{10}
\providecommand{\url}[1]{\texttt{#1}}
\providecommand{\urlprefix}{URL }
\expandafter\ifx\csname urlstyle\endcsname\relax
  \providecommand{\doi}[1]{doi:\discretionary{}{}{}#1}\else
  \providecommand{\doi}{doi:\discretionary{}{}{}\begingroup
  \urlstyle{rm}\Url}\fi

\bibitem{Bar-10}
\emph{A.~Bartel}, Large {S}elmer groups over number fields, Math. Proc.
  Cambridge Philos. Soc., \textbf{148(01)}, (2010), 73--86.

\bibitem{BB-04}
\emph{W.~Bley and R.~Boltje}, Cohomological mackey functors in number theory,
  J. of Number Theory, \textbf{105}, (2004), 1--37.

\bibitem{Bol-97}
\emph{R.~Boltje}, Class group relations from {B}urnside ring idempotents, J. of
  Number Theory, \textbf{66}, (1997), 291--305.

\bibitem{Bra-51}
\emph{R.~Brauer}, Beziehungen zwischen {K}lassenzahlen von {T}eilk\"{o}rpern
  eines {G}aloisschen {K}\"{o}rpers, Math. Nachr., \textbf{4}, (1951),
  158--174.

\bibitem{CN}
\emph{L.~Caputo and F.~Nuccio}, On fake $\mathbb{Z}_p$ extensions of number
  fields, arXiv:0807.1135v2 [math.NT].

\bibitem{CR}
\emph{C.~Curtis and I.~Reiner}, \emph{Methods of Representation Theory with
  Applications to Finite Groups and Orders}, volume~II.
\newblock John Wiley and Sons (1987).

\bibitem{Smi-01}
\emph{B.~de~Smit}, Brauer-{K}uroda relations for {S}-class numbers, Acta
  Arith., \textbf{98(2)}, (2001), 133--146.

\bibitem{Dir-42}
\emph{G.~L. Dirichlet}, Recherches sur les formes quadratiques \`{a}
  co\"{e}fficients et \`{a} ind\'{e}termin\'{e}es complexes, J. Reine Angew.
  Math., \textbf{24}, (1842), 291--371.

\bibitem{TVD-2}
\emph{T.~Dokchitser and V.~Dokchitser}, Regulator constants and the parity
  conjecture, Invent. Math., \textbf{178(1)}, (2009), 23--71.

\bibitem{TVD-1}
\emph{T.~Dokchitser and V.~Dokchitser}, On the {B}irch-{S}winnerton-{D}yer
  quotients modulo squares, Annals of Math., \textbf{172(1)}, (2010), 567--596.

\bibitem{HK-77}
\emph{F.~Halter-Koch}, Einheiten und {D}ivisorenklassen in {G}aloisschen
  algebraischen {Z}ahlk\"{o}rpern mit {D}iedergruppe der {O}rdnung 2l f\"{u}r
  eine ungerade {P}rimzahl l, Acta Arith., \textbf{33}, (1977), 353--364.

\bibitem{Kur-50}
\emph{S.~Kuroda}, \"{U}ber die {K}lassenzahlen algebraischer {Z}ahlk\"{o}rper,
  Nagoya Math. J., \textbf{1}, (1950), 1--10.

\bibitem{Lee-64}
\emph{M.~P. Lee}, Integral representations of dihedral groups of order 2$p$,
  Trans. American Math. Soc., \textbf{110(2)}, (1964), 213--231.

\bibitem{Lem-94}
\emph{F.~Lemmermeyer}, Kuroda's class number formula, Acta Arith., \textbf{66},
  (1994), 245--260.

\bibitem{Lem-05}
\emph{F.~Lemmermeyer}, Class groups of dihedral extensions, Math. Nachr.,
  \textbf{278}, (2005), 679--691.

\bibitem{Mos-75}
\emph{N.~Moser}, Unit\'{e}s et nombre de classes d'une extension galoisienne
  di\'{e}drale de $\mathbb{Q}$, Ast\'{e}risque, \textbf{24/25}, (1975), 29--35.

\bibitem{Ser-67}
\emph{J.-P. Serre}, \emph{Repr{\'e}sentation Lin{\'e}aires des Groupes finis}.
\newblock Herman Paris (1967).

\bibitem{Tat-84}
\emph{J.~Tate}, \emph{Les Conjectures de Stark sur les Fonctions L d'Artin en
  s=0}, volume~47 of \emph{Progress in Mathematics}.
\newblock Birkh\"{a}user (1984).

\end{thebibliography}
\textsc{Department of Mathematics, POSTECH
\newline San 31, Hyojadong, Namgu, Pohang, Gyungbuk 790-784
\newline Republic of Korea}
\newline E-mail address: \texttt{bartel@postech.ac.kr}
\end{document}